\documentclass[10pt]{amsart}
\allowdisplaybreaks[4]
\usepackage{amssymb,color}
\usepackage{amsthm}

\definecolor{c20}{rgb}{0.,0.7,0.}
\definecolor{c30}{rgb}{0.,0.,1.}
\definecolor{c40}{rgb}{1,0.1,0.7}
\definecolor{c50}{rgb}{1,0,0}
\definecolor{c60}{rgb}{0,0.9,0.1}

\newcommand{\E}[1]{\mathbb{E}\left\{ #1\right\}}

\newcommand{\pk}[1]{\mathbb{P} \left(#1 \right) }
\newcommand{\pl}[1]{\mathbb{P}\left(#1 \right)}

\newcommand{\BQN}{\begin{eqnarray}}
\newcommand{\EQN}{\end{eqnarray}}
\newcommand{\BQNY}{\begin{eqnarray*}}
\newcommand{\EQNY}{\end{eqnarray*}}

\newcommand{\BS}{\begin{sat}}
\newcommand{\ES}{\end{sat}}
\newcommand{\BT}{\begin{theo}}
\newcommand{\ET}{\end{theo}}
\newcommand{\BK}{\begin{korr}}
\newcommand{\EK}{\end{korr}}

\newcommand{\BD}{\begin{de}}
\newcommand{\ED}{\end{de}}

\newcommand{\BIT}{\begin{itemize}}
\newcommand{\EIT}{\end{itemize}}
\newcommand{\BDI}{\begin{description}}
\newcommand{\EDI}{\end{description}}

\newcommand{\BRM}{\begin{remarks}}
\newcommand{\ERM}{\end{remarks}}

\newcommand{\BEL}{\begin{lem}}
\newcommand{\EEL}{\end{lem}}

\newtheorem{theo}{Theorem}[section]
\newtheorem{sat}[theo]{Proposition}
\newtheorem{de}[theo]{Definition}
\newtheorem{lem}[theo]{Lemma}

\newtheorem{korr}[theo]{Corollary}
\newtheorem{remark}[theo]{Remark}
\newtheorem{remarks}[theo]{Remarks}

\newcommand{\prooflem}[1]{\textbf{Proof of Lemma} \ref{#1}}

\newcommand{\COM}[1]{}

\newcommand{\QED}{\hfill $\Box$ \\}

\def\rw{\rightarrow}

\def\IF{\infty}

\def\H{\mathcal{H}}

\def\kd#1{\textcolor{black}{#1}}
\def\rd#1{\textcolor{black}{#1}}
\def\gr#1{\textcolor{black}{#1}}
\def\rdd#1{\textcolor{black}{#1}}
\def\bl#1{\textcolor{black}{#1}}

\def\rd#1{\textcolor{black}{#1}}

\begin{document}

\title
{  Extremes of nonstationary  Gaussian fluid queues}

\author{Krzysztof D\c{e}bicki}
\address{Krzysztof D\c{e}bicki, Mathematical Institute, University of Wroc\l aw, pl. Grunwaldzki 2/4, 50-384 Wroc\l aw, Poland}
\email{Krzysztof.Debicki@math.uni.wroc.pl}
\author{Peng Liu}
\address{Peng Liu, Department of Actuarial Science, University of Lausanne, UNIL-Dorigny 1015 Lausanne, Switzerland}
\email{peng.liu@unil.ch}

\bigskip

%\date{\today}
 \maketitle
\bigskip
{\bf Abstract:}
This contribution investigates asymptotic properties of transient
queue length process
\[
Q(t)=\max\left(Q(0)+X(t)-ct, \sup_{0\leq s\leq t}\left(X(t)-X(s)-c(t-s)\right)\right),\ \ \ t\geq 0
\]
in Gaussian fluid queueing model, where
input process $X$ is modeled by a centered Gaussian process with stationary increments and
$c>0$ is the output rate.
More specifically, under some mild conditions on $X$ and $Q(0)=x\ge0$, the exact asymptotics
of
$$\pi_{x,T_u}(u)=\pk{Q(T_u)>u},$$
as $u\to\infty$, is derived.
The play between $u$ and $T_u$ leads to two qualitatively different regimes:
(A) short-time horizon when $T_u$ is relatively small with respect to $u$;
(B) moderate- or long-time horizon when $T_u$ is asymptotically much larger than $u$.
As a by-product, some implications for the speed of convergence to stationarity of the considered model
are discussed.



\section{Introduction}
The analysis of queueing systems with Gaussian input attracted substantial
interest in last years. The importance of modelling input stream by a Gaussian process stems both from
theory-oriented arguments, mainly based on {\it central limit theorem}-type results
applied to multiplexed input streams (see, e.g., \cite{ DE1999, MiS07, TWS1997, WTSW1995}) and applied-oriented approach taking
advantage of richness and flexibility of the class of Gaussian processes, allowing to model such phenomena as
{\it long range dependence} or {\it self-similarity}.

Consider queue fed by a Gaussian process with stationary increments $X(t)$
and emptied at rate $c>\E {X(1)}$.  Having the interpretation that, for $s<t$,
$X(t)-X(s)$ is the amount of traffic having entered to the system in time interval $[s,t)$,
we define the buffer content process $\{Q(t),t\ge0\}$
by %the following representation
\BQN
Q(t)=\max\left(Q(0)+X(t)-ct, \sup_{0\leq s\leq t}\left(X(t)-X(s)-c(t-s)\right)\right),\ \ \ t\geq 0.
\label{q.q0}
\EQN

Vast majority of literature on properties of $Q(t)$ deals with the %corresponding
{\it steady-state} solution of (\ref{q.q0}), which takes form
\BQN
Q^*(t)=\sup_{-\IF \le s\le t}(X(t)-X(s)-c(t-s)), \ \ t\geq 0, \label{q.stat}
\EQN
%see, e.g., \cite{Pit2001}.
with particular focus on the asymptotics of probability
that the steady state buffer content exceeds high level $u$, that is
\BQN
\widehat{\pi}(u):= \pl{Q^*(0)>u}\label{st.1},\ \  u\to\infty,
\EQN
see \cite{Norros94,HP1999,DEK2002,HP2004,DI2005} and references therein.
We refer also to  counterparts of (\ref{st.1}) under {\it many-source} regime
(e.g. \cite{DeMan03} or
monograph \cite{Man07})
and related
recent results on asymptotics for extremes of $\gamma-$reflected Gaussian
processes \cite{KEP2015,HA2013}.
\\

Substantially less is known on
nonstationary characteristics of the queue content process (\ref{q.q0}), in particular
if $Q(0)>0$.
%,  transient counterparts of
%(\ref{st.1}) and (\ref{st.2}).}
%The reason for this is that in
In this case the system
additionally depends on
the initial queue content at time $t=0$ and
on time at which it is analyzed, leading to more complicated structure of
the queue process, which makes the analysis of the distribution of (\ref{q.q0})
much more difficult.
%We refer to \cite{Deb03} for the analysis of a special case
%with
%The topic of this paper stems straightforwardly from this open questions.
%More specifically, let us
More specifically, suppose that $Q(0)=x\ge0$ and rewrite (\ref{q.q0}) as
\BQN
Q(t)=\max\left(x+X(t)-ct, \sup_{0\leq s\leq t}\left(X(t)-X(s)-c(t-s)\right)\right),\ \ \ t\geq 0.
\label{main}
\EQN
This contribution is devoted to the analysis of the exact asymptotics
of the tail distribution of
the nonstationary workload $Q(t)$ defined by (\ref{main}) at time $T_u$, i.e.
\begin{eqnarray}
\pi_{x,T_u}(u):=\pl{Q(T_u)>u}, \label{main.1} \ {\rm as} \ u\to\infty.
\end{eqnarray}
It appears that the play between $x, T_u$ and $X$ leads
to several scenarios which can be grouped according to the relation between
$T_u$ and $u$ on: {\it short-time} and {\it moderate- or long-time} case.
Then, within each of the above time-horizons, one faces several types of the asymptotics.
The \rd{ results derived  in this contribution} complement
findings obtained for the stationary systems,
see \cite{Norros94,HP1999,HP2004,DEK2002,DI2005,KrzysPeng2015}
and
extend results of \cite{DebickiRol02}, where (\ref{main.1})
was considered for much simpler case
$x=0$ and $T_u=T>0$.
%We note that
The complexity of the derivations of main results of
this contribution is substantially higher than that of the corresponding proofs in the above papers.
More specifically, \rd{with $\Psi(\cdot)$ being  the tail distribution
function of a standard Gaussian random variable,} using that
\BQN\label{Basic1}
\pi_{x,T_u}(u)&=&\pi_{0,T_u}(u)+\Psi\left(\frac{u-x+cT_u}{\sigma(T_u)}\right)\\
&& \ \ -\pl{X(T_u)-c T_u>u-x, \sup_{0\leq s\leq T_u}\left(X(T_u)-X(s)-c(T_u-s)\right)>u}\nonumber,
\EQN
 one can distinguish two scenarios.
First, if $\pi_{0,T_u}(u)$ or $\Psi\left(\frac{u-x+cT_u}{\sigma(T_u)}\right)$
is asymptotically dominating, then by (\ref{Basic1}) it determines the asymptotics
of $\pi_{x,T_u}(u)$ as $u\rw\IF$.
Utilizing that
$ \pi_{0,T_u}(u)=\pl{\sup_{t\in[0,T_u]}X(t)-ct>u}$,
the main idea of the proofs in this case is based on an extension of the {\it double sum} method, a technique
which was originally developed for the study of
asymptotics of suprema of centered Gaussian processes; see e.g. \cite{PicandsB, PicandsA, MR0494458}
and
monographs \cite{Pit96, Pit20}.
Second, when
$\pi_{0,T_u}(u)$ is asymptotically of the same order as $\Psi\left(\frac{u-x+cT_u}{\sigma(T_u)}\right)$, as $u\to\infty$,
\rd{then one}
needs an independent approach that goes beyond the
double sum method and
leads to new types of asymptotics that are not present in the literature on Gaussian extremes
(see Section \ref{s.x0}).

The model analyzed in this paper covers
wide class of Gaussian inputs, including the celebrated {\it fractional Brownian
motions} and {\it Gaussian integrated} processes.

The derived results shed some light on important
issues related to the speed of convergence to stationarity
of the queueing system in time; see \cite{MR2538075}
for works with fractional Brownian motion input.
In particular, by comparing our findings with their counterparts
for the stationary model,
we arrive at \rd{a  finding}
that
%$\diamond$ For
the system which starts off with empty queue asymptotically
(for large $u$)
 \rd{reaches the steady state asymptotics faster than} the nonempty one.

A related problem that addresses transient properties of the buffer content process
is the analysis of $Q(T)$ conditioned by its initial content $Q(0)$.
We refer to \cite{DSM09}, where the logarithmic asymptotics of
$\pk{ Q^*(T)>pu, Q^*(0)>qu}$ as $u\to\infty$, for $p,q>0$  was derived,
giving some insight into the asymptotics of $\pk{ Q^*(T)>pu |Q^*(0)>qu}$
for the stationary buffer content process $Q^*$.

Organisation of the paper: Section \ref{s.model} contains introduction of the model and
the notation.
The main results of the paper are presented in Section \ref{s.x0}.
Section \ref{s.conv} is devoted to
relation between
(\ref{main.1}) and
(\ref{st.1}).
Some technical results that are useful in the proofs are given in
Section \ref{s.aux}.
Section \ref{s.proofs} contains
detailed proofs of the  main results.

%In this contribution we analyze the behavior of the system
%which didn't reach stationary equilibrium, focusing on
%the asymptotic behavior of the queueing process $Q_x(t)$, where $x\ge0$
%is the buffer content at time $t=0$.

\section{Model description and preliminary results}
\label{s.model}

Let $\{X(t),t\geq 0\}$ be a Gaussian process with \rd{stationary increments, $X(0)=0$ a.s. and variance function } $\sigma^2(t)$.
With no loss of generality we assume \kd{that $X(t)$ is} centered, i.e.
$\E {X(t)}\equiv 0$, $t\ge0$.
We suppose that\\ %Additionally, its variance function satisfies\\
\\
{\bf AI}:
$\sigma^2(t)\sim A_\IF t^{2\alpha_\IF}$ as $t\rw\IF$  with $A_\IF>0, \alpha_\IF\in (0,1)$. Further, $\sigma^2(t)$ is twice continuously differentiable on $(0,\IF)$ with its first derivative $\dot{\sigma^2}$ and second derivative $\ddot{\sigma^2}$ being ultimately monotone at $\IF$.\\
{\bf AII} $\sigma^2(t)\sim A_0 t^{2\alpha_0}$ as $t\rw 0$ with $A_0>0, \alpha_0\in(0,1]$.\\

We note that assumptions {\bf AI-AII} cover \kd{all classical Gaussian} input models considered
in the literature, including {\it fractional Brownian motion} $X(t)=B_H(t)$
(i.e. $Var(X(t))=t^{2H}$, with $H\in(0,1)$),
see \cite{Norros94,HP1999}, and {\it integrated Gaussian} inputs,
where $X(t)=\int_0^t Z(s)ds$, with $Z(t)$ being a centered stationary Gaussian process with covariance function
satisfying some standard conditions; see e.g. \cite{DEK2002,HP2004,DI2005}.
%Suppose that
%\BQN\label{SUP}
%\CE{X(t)}=\int_0^tY(s)ds, t\geq 0,
%\EQN
% where $Y$ is a stationary
% \KD{centered}
% Gaussian process with continuous trajectories.
%Hence the ruin probability is
%$$\psi_{\gamma,T}\eL{(u)}=\mathbb{P}\left(\sup_{0\leq s\leq t\leq T}\left(Y_I(t)-\gamma Y_I(s)-c(t-\gamma s)\right)>u\right)
%$$
%for $T\in(0,\IF]$.
%Let $R(t)$ denote the correlation function of $Y$ and without loss of generality, we suppose that $R(0)=1$. In this subsection, we shall consider two scenarios:\\
% {\bf SRD} (short-range dependent), i.e., we shall assume that\\
% i) $R(t)$ is ultimately monotone and $\lim_{t\rw\IF}t r(t)=0$,\\
% ii) %$\int_0^\IF|R(t)|dt<\IF$ and
% \KD{$\int_0^\IF R(t)dt=G \in \ehhe{(0,\IF)}$}.\\
% {\bf LRD} (long-range dependent), {i.e., we shall suppose that} \\
% i)   $R(t)$ is decreasing over $[0,\IF)$,\\
% ii) $R(t)\sim \vartheta t^{2H-2}$ as $t\rw\IF$ with $\vartheta>0, H\in (1/2, 1)$.\\
%It follows that {\bf AI-AII} are satisfied if  $X$ is {\bf SRD} or {\bf LRD}.

Following the introduction, we consider
a queue fed by input process $X(t)$ and emptied at a constant rate $c>0$.
The queue content process $Q(t)$, with $x=Q(0)\ge0$,
is defined as in (\ref{main}).
%is the initial content of the system at time $t=0$
%by
%\BQN
%Q_x(t):=\max\left(x+X(t)-ct, \sup_{0\leq s\leq t}\left(X(t)-X(s)-c(t-s)\right)\right),\ \ \ t\geq 0.
%\label{main}
%   \EQN

Due to (\ref{Basic1}),
%\BQN\label{Basic1}
%\pi_{x,T_u}(u)&=&\pi_{0,T_u}(u)+\Psi\left(\frac{u-x+cT_u}{\sigma(T_u)}\right)\\
%&& \ \ -\pl{X(T_u)-c T_u>u-x, \sup_{0\leq s\leq T_u}\left(X(T_u)-X(s)-c(T_u-s)\right)>u}\nonumber,
%\EQN
for the analysis of (\ref{main.1})
it is convenient %to focus first on the case $x=0$.
%Thus in Section \ref{s.x0} we
to start with detailed asymptotic analysis of  $\pi_{0,T_u}(u)$
 as $ u\to\infty$.
% distinguishing relation between $T_u$ and $u$.
%The results obtained in this section complement findings of \cite{DebickiRol02} and \cite{GlZ},
%by giving exact aymptotics and allowing the time horizon to depend on the threshold $u$.

Having that
\[
\pi_{0,T_u}(u)=\pl{\sup_{t\in[0,T_u]} X(t)-ct>u}=
\pl{\sup_{t\in[0,T_u]} \frac{X(t)}{u+ct}>1},
\]
function
\BQN\label{mut}
m(u,t):=Var^{-1/2}\left(  \frac{X(t)}{u+ct} \right)= \frac{u+ct}{\sigma(t)}
\EQN
will play crucial role in further analysis.
Let
\BQN\label{tu}\kd{t_u}:=\arg \min_{t \ge 0} m(u,t)\EQN
and observe that under {\bf AI-AII}, as shown in Lemma \ref{L2}, we have
\BQN\label{t*}
\kd{t^*}:=\lim_{u\to\infty} \frac{\kd{t_u}}{u}=\frac{\alpha_\infty}{c(1-\alpha_\infty)}.
\EQN
Next, we introduce
\BQN\label{delta}
\Delta(u,s)=
\overleftarrow{\sigma}
 \left(\frac{\sqrt{2}\sigma^2(s)}{u+cs}\right)
 %\sim\left\{\begin{array}{cc}
%\left(\frac{\sqrt{2A_\IF}(\tau^*)^{2\alpha_\IF}}{1+c\tau^*}\right)^{1/\alpha_\IF} u^{\frac{2\alpha_\IF-1}{\alpha_\IF}}, & \alpha_\IF>1/2\\
%\overleftarrow{\sigma}
% \left(\frac{\sqrt{2}A_\IF\tau^*}{1+c\tau^*}\right), & \alpha_\IF=1/2\\
% \left(\frac{\sqrt{2}A_\IF(\tau^*)^{2\alpha_\IF}}{\sqrt{A_0}(1+c\tau^*)}\right)^{1/\alpha_0} u^{\frac{2\alpha_\IF-1}{\alpha_0}}, & \alpha_\IF<1/2\\
% \end{array}\right.
\EQN
with $\overleftarrow{\sigma}$ the asymptotic inverse function of $\sigma$.
Moreover, we \rd{denote}
\BQN\label{AB}
A=\left(\frac{{\alpha_\IF}}{c(1-{\alpha_\IF})}\right)^{-\alpha_\IF}\frac{1}{1-\alpha_\IF}, \quad
B=\left(\frac{\alpha_\IF}{c(1-\alpha_\IF)}\right)^{-\alpha_\IF-2}\alpha_\IF.
\EQN

Finally, we introduce constants that appear in the derived asymptotics.
Let
\BQNY
\mathcal{H}_{X}[0,S]=\E{e^{\sup_{t\in[0, S]}\left(
 \sqrt{2} X(t)-\sigma^2(t)\right)}},\ \ \mathcal{P}_{X}^{f}[0,S]=\E{e^{\sup_{t\in[0, S]}\left(
 \sqrt{2} X(t)-\sigma^2(t)-f(t) \right)}},
 \EQNY
where $X$ is a centered Gaussian process with stationary increments that satisfies {\bf AI-AII}, $f$ is a nonnegative function over $[0,\IF)$ and  $S>0$.
The {\it generalized Pickands} and {\it Piterbarg constants} are defined by
 \BQNY
 \mathcal{H}_X:=\lim_{S\rw\IF}\frac{\mathcal{H}_{X}[0,S]}{S}, \ \ \mathcal{P}_{X}^{f}:=\lim_{S\rw\IF} \rd{\mathcal{P}_{X}^{f}[0,S]}
 \EQNY
 respectively.
 We refer to, e.g.,  \cite{Michna1, KEP2015, DEK2002,  Harper2,DM, DiekerY, DE2014,Harper3, En16, Michna2, Shao}
for the proof of existence and properties of
(generalized) Pickands and Piterbarg constants, simulation issues and their relations to
max-stable processes.
Additionally, %the following new {\it Pickands-Piterbarg} type constants which will appear in the present paper.
%for ${X}$, an independent copy of Gaussian process $X$,
%$X^*(t):=\sqrt{2}X(t)-\sigma_{X}^2(t)-f(t)$ and $\widetilde{X}^*(t):=\sqrt{2}\widetilde{X}(t)-\sigma_{\widetilde{X}}^2(t)-f(t)$, we introduce
for a given nonnegative function $f$ and $a\ge0$,
let
$$\mathcal{P}_{X}^{f,a}[0,S]=\E{e^{\max\left(a, \sup_{t\in[0,S]}(\sqrt{2}X(t)-\sigma_{X}^2(t)-f(t))\right)}},$$
%$$\widetilde{\mathcal{P}}_{X}^{f, a}[0,S]=\E{e^{\max\left(\sup_{t\in[0,S]}a+X^*(t), \sup_{s,t\in[0,S]}X^*(t)+{\widetilde{X}}^*(s)\right)}},$$
and
$$\mathcal{P}_{X}^{f, a}:=\lim_{S\rw\IF}\mathcal{P}_{X}^{f, a}[0,S].$$
These constants appear in the asymptotics of
(\ref{main.1}) for some   \rd{scenarios considered in the next section.}
\rd{Note that, for any $S>0$,
\BQNY
\mathcal{P}_{X}^{f, a}[0,S]\leq e^a \mathcal{P}_{X}^{f}[0,S].
\EQNY
This implies that if $\mathcal{P}_{X}^{f}<\IF$, then
$$\mathcal{P}_{X}^{f, a}\leq e^a \mathcal{P}_{X}^{f}<\IF.$$ }
Let $\Phi(\cdot)$  be the distribution
function of a standard Gaussian random variable and
we write $f(u)\sim g(u)$ to denote the asymptotic equivalence $\lim_{u\to\infty}\frac{f(u)}{g(u)}=1$.
Before proceeding to main results of this paper that deal with the case where $T_u\to\infty$, as $u\to\infty$,
we provide a preliminary \rd{one} that covers the easier case $T_u=T>0$.

\BS \label{p.T}
Suppose that
$\dot{\sigma^2}(t)>0$ for $t\in (0,T]$ and {\bf AII} holds. Then for $T\in (0,\IF)$ and $x>0$, as $u\rw\IF$,
$$\pi_{x,T}(u)\sim \Psi\left(\frac{u-x+cT}{\sigma(T)}\right).$$
\ES
The proof of Proposition \ref{p.T} is deferred to Section \ref{s.proofs}.
In next section we tacitly assume that $T_u\to\infty$, as $u\to\infty$.

\section{Main results}
\label{s.x0}
\kd{The asymptotics of $\pi_{x, T_u}(u)$, as $u\to\infty$,}
strongly depends
on the relation between $T_u$ and $u$, leading to two separate scenarios: i) short-time horizon
and ii) moderate- or long-time horizon, which we analyze separately.

\subsection{Short-time horizon}\label{s.short}
In this section we consider the case \kd{where} $T_u$ is relatively small with comparison to $u$.
More precisely, we suppose that $T_u\to\infty$, as $u\to\infty$, and\\
%\begin{eqnarray}
{\bf T1} $\lim_{u\to\infty}\frac{T_u}{u}=\gamma \in [0,t^*)$.%\label{short}
\\
%\end{eqnarray}
%Using that under {\bf T1}, $T_u<t^* u \sim t_u$ for large $u$, intuitively
%$T_u$ is not sufficiently large for the system to reach the steady state in this scenario.
Let
\BQN\label{varphi}
\varphi:=\lim_{u\to\infty}\frac{T_u}{u^{1/(2\alpha_\infty)}}\in[0,\infty]
\EQN
and %impacts the form of the derived asymptotics. %Define
define \rd{a family of Gaussian random processes $\{\mu_\varphi(t), t\in\mathbb{R}\}$, where }
\BQN\label{muvarphi}\mu_\varphi(t)=\left\{\begin{array}{cc}
B_{\alpha_0}(t), & \hbox{if } \varphi=0\\
 \frac{1+c{\gamma}}{\sqrt{2}A_\IF\varphi^{2\alpha_\IF}}X
\left(\overleftarrow{\sigma}\left(\frac{\sqrt{2}A_\IF\varphi^{2\alpha_\IF}}{1+c{\gamma}}\right)t\right), & \hbox{if } \varphi\in (0,\IF)\\
B_{\alpha_\IF}(t), & \hbox{if } \varphi=\IF.
\end{array}\right.
\EQN

%$$\mu_{\varphi}(t)=\left\{
 %           \begin{array}{ll}
%B_{\alpha_0}(t), & \hbox{if } \varphi=0 \\
%\frac{1+c{\gamma}}{\sqrt{2}A_\IF\varphi^{2\alpha_\IF}}X\left(\overleftarrow{\sigma}(\frac{\sqrt{2}A_\IF \varphi^{2\alpha_\IF}}{1+c\gamma})t\right), & \hbox{if }\varphi\in (0,\IF)\\
%B_{\alpha_\IF}(t),&  \hbox{if } \varphi=\IF.
%              \end{array}
 %           \right.$$
%\subsubsection
Due to equation (\ref{Basic1})
it is convenient first to analyze the system that starts off with empty queue.
\\
\\
$\diamond$ {\underline{{Case $Q(0)=0$.}}
In this scenario
the asymptotic behaviour of $\Omega(u,T_u)$
as $u\to\infty$, where
\BQN\label{omega}
\Omega(u,t):=\frac{m^2(u,t)\Delta(u,t)}{t}, \quad t>0,
\EQN
 leads to three qualitatively different cases.

\BT\label{submain1}
Suppose that $T_u$ satisfies {\bf T1}.\\
%=\gamma_\beta u^\beta(1+o(1))$ with $\gamma\in (0,1]$ and $\gamma_\beta$ satisfying (\ref{gamma}).\\
i) If $\lim_{u\rw\IF}\Omega(u,T_u)=0$, then
\BQNY
\pi_{0,T_u}(u)&\sim&\mathcal{H}_{B_{\alpha_0}}\left(\alpha_\IF-\frac{c{\gamma}}{1+c{\gamma}}\right)^{-1}
%\frac{T_u}{\Delta(u,T_u)m^2(u,T_u)}
\Omega^{-1}(u,T_u)
\Psi(m(u,T_u)).
\EQNY
ii) If $\lim_{u\rw\IF}\Omega(u,T_u)=\Omega_\IF \in(0,\IF)$, then
\BQNY
\pi_{0,T_u}(u)
\sim \mathcal{P}_{\mu_\varphi}^{\Omega_\IF(\alpha_\IF-\frac{c{\gamma}}{1+c{\gamma}})t}\Psi(m(u,T_u))
%\left\{\begin{array}{ll}
%\mathcal{P}_{B_{\alpha_0}}^{\Omega_\IF(\alpha_\IF-\frac{c{\gamma}}{1+c{\gamma}})t}\Psi(m(u,T_u)), & \hbox{if }\varphi=0\\
%\mathcal{P}_{\mu}^{\Omega_\IF(\alpha_\IF-\frac{c{\gamma}}{1+c{\gamma}})t}\Psi(m(u,T_u)), & \hbox{if }\varphi\in (0,\IF)
%\end{array}\right.
\EQNY
with
$\varphi\in [0,\IF)$.\\
iii) If $\lim_{u\rw\IF}\Omega(u,T_u)=\IF$, then
\BQNY
\pi_{0,T_u}(u)\sim\Psi(m(u,T_u)).
\EQNY
\ET
%As an immediate conclusion of the above theorem, we arrive at
%\BK\label{Cor1} If  {\bf T1} is satisfied,  then
%\BQNY
%\pi_{0,T_u}(u-x)\sim\pi_{0,T_u}(u)e^{\frac{u+cT_u}{\sigma^2(T_u)}x}.
%\EQNY
%\EK

%$
%\mu(t)=\frac{1+c{\gamma}}{\sqrt{2}A_\IF\varphi^{2\alpha_\IF}}X
%\left(\overleftarrow{\sigma}(\frac{\sqrt{2}A_\IF\varphi^{2\alpha_\IF}}{1+c{\gamma}})t\right)
%$.\\

\vspace{0.5cm}

$\diamond$ {\underline{{Case $Q(0)>0$.}}
Now we analyze asymptotic properties of the system that starts off with nonempty queue;
up to the end of this subsection we tacitly suppose that $x=Q_x(0)>0$.
In order to make the results of this case more transparent,
we present the derived asymptotics in
the language of $\pi_{0,T_u}(u)$ and $\Psi\left(\frac{u-x+cT_u}{\sigma(T_u)}\right)$ respectively, which was
derived in previous section.

%\vspace{0.2cm}
We begin with
derivation of the asymptotics of
$\pi_{x,T_u}(u)$, as $u\to\infty$.
Following observation (\ref{Basic1}),
if one of terms
$\Psi\left(\frac{u-x+cT_u}{\sigma(T_u)}\right)$ or
$\pi_{0,T_u}(u)$ is asymptotically dominant, then it rules the asymptotic behavior of
$\pi_{x,T_u}(u)$.
The situation when $\Psi\left(\frac{u-x+cT_u}{\sigma(T_u)}\right)=O(\pi_{0,T_u}(u))$
is particularly delicate. Its proof needs a case-specific analysis and
leads to a separate form of the asymptotics;
see case $\varphi\in (0,\infty)$
%and $\alpha_\infty=1/2$
in theorem below.
%\K{Let, for $\varphi\in (0,\IF)$, .}
We note that assumption {\bf T1} together with $\varphi\in(0,\infty)$
implies that $\alpha_\IF\ge1/2$.

\BT\label{nth1}
Suppose that $T_u$ satisfies {\bf T1}.\\
i) If $\varphi=0$, then
\[
\pi_{x,T_u}(u)\sim \Psi\left(\frac{u-x+cT_u}{\sigma(T_u)}\right).
\]
ii)If $\varphi\in(0,\infty)$, then
\begin{eqnarray*}
\pi_{x,T_u}(u)
\sim
\left\{\begin{array}{ll}
\Psi\left(\frac{u-x+cT_u}{\sigma(T_u)}\right)
%\sim e^{\gamma_{1/(2\alpha_\IF)}^{-2\alpha_\IF}\frac{x}{A_\IF}}\pi_{0,T_u}^{\sup}(u)
, & \hbox{if }\alpha_\IF>1/2\\
\mathcal{P}_{a_1X}^{a_2t, \sqrt{2}a_1x}\Psi(m(u,T_u)), & \hbox{if }\alpha_\IF=1/2,
%\sim
%\widetilde{\mathcal{P}}_{a_1X}^{a_2t, \sqrt{2}a_1x}\left({P}_{a_1X}^{a_2t}\right)^{-2}\pi_{0,T_u}^{\sup}(u)\sim
%\widetilde{\mathcal{P}}_{a_1X}^{a_2t, \sqrt{2}a_1x}\left({P}_{a_1X}^{a_2t}\right)^{-1}e^{-\sqrt{2}a_1}\pi_{0,T_u}(u-x)
\end{array}\right.
\end{eqnarray*}
where
\BQN\label{a12}
a_1=\frac{1+c\varphi}{\sqrt{2}A_\IF\varphi},\quad
a_2=\frac{(1+c\varphi)^2}{A_\IF\varphi^2}\left(\alpha_\IF-\frac{c\varphi}{1+c\varphi}\right).
\EQN

iii) If $\varphi=\infty$, then
\[
\pi_{x,T_u}(u)\sim \pi_{0,T_u}(u).
\]
\ET
\vspace{0.5cm}

We observe that, under {\bf T1} combined with $\varphi=\infty$, by
Theorem \ref{nth1}, the asymptotics of $\pi_{x,T_u}(u)$, as $u\to\infty$,
doesn't depend on the initial buffer content $x$.

\subsection{Moderate- and long-time horizon}
Now, let us proceed to the case that $T_u$ is "moderate" or "large"
with comparison to $u$. To be more precise, in this section we suppose that
%\begin{eqnarray}
\\
\\
{\bf T2} $\lim_{u\to\infty} \frac{T_u-t_u}{u^{\alpha_\infty}} = \omega\in(-\infty,\infty]$.
\\
\\
Recall that $t_u\sim t^*u=\frac{\alpha_\infty}{c(1-\alpha_\infty)}u$, as $u\to\infty$.
We note that, if $\omega\in (-\infty,\infty)$, then $T_u$ is asymptotically close to $t_u$ (moderate-time horizon), while
$\omega=\infty$ deals with the case where $T_u$ is relatively large with comparison to $\kd{t_u}$ (long-time horizon).

Let
\BQN\label{eta}
\eta_{\alpha_\IF}(t)=\left\{
            \begin{array}{ll}
B_{\alpha_0}(t), & \hbox{if } \alpha_\IF<1/2 \\
\frac{1+ct^*}{\sqrt{2}A_\IF t^*}X\left(\overleftarrow{\sigma}\left(\frac{\sqrt{2}A_\IF t^*}{1+ct^*}\right)t\right), & \hbox{if }\alpha_\IF=1/2\\
B_{\alpha_\IF}(t),&  \hbox{if } \alpha_\IF>1/2.
              \end{array}
            \right.
\EQN

Analogously to the short-time horizon scenario investigated in Section
\ref{s.short}, we separately consider the case of empty and nonempty system at $t=0$.
\\
\\
$\diamond$ {\underline{{Case $Q(0)=0$.}}
We begin with the asymptotic analysis of
$\pi_{0,T_u}(u)$ under {\bf T2}.% takes uniform form, although depending on $\omega$.

\BT\label{th.0.2}
Suppose that
$T_u$ satisfies {\bf T2}.
Then
\BQNY
\pi_{0,T_u}(u)&\sim& \mathcal{H}_{\eta_{\alpha_\IF}}\sqrt{\frac{2A\pi}{B}}
\frac{u}{m(u,t_u)\Delta(u,t_u)}
\Phi\left(\sqrt{\frac{B}{AA_\IF}}\frac{(1+ct^*)\omega}{(t^*)^{\alpha_\IF}}\right)
\Psi(m(u,t_u)).
\EQNY
\ET

%It is interesting to compare the above asymptotics with the asymptotics of
%the overflow probability in the stationary system....
%\BK

\vspace{0.5cm}
$\diamond$ {\underline{{Case $Q(0)>0$.}}
Up to the end of this section we suppose that the queue is nonempty at time $t=0$, i.e.
$x=Q_x(0)>0$.
It appears that, under {\bf T2}, this scenario
delivers qualitatively different types of the asymptotics than the case $Q(0)=0$.
\BT\label{nth3}
Suppose that
$T_u$ satisfies {\bf T2}.\\
i) If $\alpha_\infty<1/2$ and
${ \lim_{u\to\infty} \frac{T_u-t_u}{\sqrt{u}}=\vartheta\in [0, \infty]}$, then
\begin{eqnarray*}
\pi_{x,T_u}(u)\sim
\left\{
            \begin{array}{ll}
\Psi\left(\frac{u-x+cT_u}{\sigma(T_u)}\right), & \hbox{if } \vartheta<\sqrt{ \frac{2A}{B}    (1-\alpha_\IF)x} ,\\
\Psi\left(\frac{u-x+cT_u}{\sigma(T_u)}\right)+ \pi_{0,T_u}(u), & \hbox{if }\vartheta=\sqrt{\frac{2A}{B}(1-\alpha_\IF)x},\\
\pi_{0,T_u}(u),&  \hbox{if } \vartheta>\sqrt{\frac{2A}{B}(1-\alpha_\IF)x} .
              \end{array}
            \right.
\end{eqnarray*}
ii)
If $\alpha_\infty\ge1/2$, then
\begin{eqnarray*}
\pi_{x,T_u}(u)
\sim
\pi_{0,T_u}(u).
\end{eqnarray*}
\ET
\vspace{0.5cm}

% ii) For the case that $T_u>0$ and $T_u\rw 0$, if $\sigma^2(t)$ is regularly varying at $0$ with index $2\alpha\in (0,2)$ or $\sigma^2(t)\sim C_2t^2$ as $t\rw 0$, then we can conclude that  $\pi_{x,T}^{\sup}(u)\sim \pi_{0,T}(u-x)$ and $\pi_{0,T}^{\sup}(u)=o(\pi_{0,T}(u-x))$. The assumptions ensure that
% \BQNY
% \E{\left(\frac{X(T_ut)-X(T_us)}{u+cT_u(t-s)}\frac{u+cT_u}{\sigma(T_u)}-\frac{X(T_ut')-X(T_us')}{u+cT_u(t'-s')}\frac{u+cT}{\sigma(T_u)}\right)^2}\leq \mathbb{Q}(|t-t'|^{\beta_4}+|s-s'|^{\beta_4}),\ \ s,t,s',t'\in[0,1]
% \EQNY
% with some $\beta_4>0$.
% Then the similar arguments as above can verify the claims.\QED

%%%%%%%%%%%%%%%%%%%%%%%%%%%%%%%%%%%%%%%%%%%%%%%%%%%%%%%%%%%%%%%%%%%%%%%%%
\section{Speed of convergence to stationarity}\label{s.conv}

This section is devoted to some remarks on
the speed of convergence of the distribution of
$Q(T_u)$ with $Q(0)=x\geq 0$ to its stationary counterpart $Q^*(0)$.
Comparison of the results derived in Section \ref{s.x0} with
asymptotics for the stationary system
given in \cite{DI2005} and \cite{KrzysPeng2015} allows us to
give some \rd{insight} into this issues.
%Following the notation introduced in Section $1$,
%let $Q^*(t)$
%denote the stationary buffer content process defined by (\ref{q.stat})
%and
Let
$
\widehat{\pi}(u):=\pk{Q^*(0)>u}.
$
Straightforward combination of results in \cite{DI2005} with Theorem \ref{th.0.2} and \gr{Theorem \ref{nth3}}
leads to the following proposition.
\BS\label{p.1}
Suppose that
$T_u$ satisfies {\bf T2}.\\
i) If $\kd{Q(0)=0}$, then
\[\pi_{0,T_u}(u)\sim \Phi\left(\sqrt{\frac{B}{AA_\IF}}\frac{(1+ct^*)\omega}{(t^*)^{\alpha_\IF}}\right) \widehat{\pi}(u).\]
ii) If $Q(0)=x>0$, $\alpha_\infty<1/2$ and $\limsup_{u\rw\IF}\frac{T_u-t_u}{\sqrt{u}}< \sqrt{   \frac{2A}{B}   (1-\alpha_\IF)x}$, then
$$\widehat{\pi}(u)=o(\pi_{x,T_u}(u)).$$
iii) If $Q(0)=x>0$, $\alpha_\infty<1/2$ and  $\liminf_{u\rw\IF}\frac{T_u-t_u}{\sqrt{u}}>\sqrt{   \frac{2A}{B}  (1-\alpha_\IF)x}$, then
$$\pi_{x,T_u}(u)\sim \widehat{\pi}(u).$$
iv) If $Q(0)=x>0$ and $\alpha_\IF\geq 1/2$, then
$$\pi_{x,T_u}(u)\sim \Phi\left(\sqrt{\frac{B}{AA_\IF}}\frac{(1+c\kd{t^{*}}
)\omega}{(t^*)^{\alpha_\IF}}\right)  \widehat{\pi}(u).$$
\ES
%As a straightforward corollary from Proposition \ref{p.1}
%we observe that in order to have $\pi_{0,T_u}(u)\sim \widehat{\pi}(u)$,
\kd{From the above proposition we see that
for $\alpha_\infty<1/2$, which corresponds to
{\it short-range dependent}
structure of the input process $X$ in the sense
that
$\sum_{k=1}^\infty {\rm Cov}(X(k)-X(k-1),X(1))<\infty$,
the system reaches stationary asymptotics faster if
it starts off with empty queue (i.e.  $Q(0)=0$) in comparison to nonempty system at time $t=0$.
For $\alpha_\infty\ge 1/2$, the initial content of the queue
doesn't influence the speed of convergence to the stationary asymptotics.}

\begin{remark}
The case
$Q(0)=x>0$, $\alpha_\infty<1/2$ and $\lim_{u\rw\IF}\frac{T_u-t_u}{\sqrt{u}}= \sqrt{   \frac{2A}{B}   (1-\alpha_\IF)x}$
is sensitive to higher order asymptotic expansion of $T_u$ and $t_u$, which needs additional knowledge on the asymptotics of
$\sigma^2(t)$. This leads to tedious calculations which go beyond the setup of this contribution.
\end{remark}

%%%%%%%%%%%%%%%%%%%%%%%%%%%%%%%%%%%%%%%%%%%%%%%%%%%%%%%%%%%%%%%%
\section{Auxiliary lemmas}\label{s.aux}
In this section we display  technical lemmas that will be helpful in the forthcoming  proofs.
\bl{ In order to improve the readability of proofs of the main results,
we list the glossary of notation that we use in the proofs.
We recall that $\sigma(t)=\sqrt{Var(X(t))}$, $c$ is given in (\ref{main}), and  $A_0, A_\IF$ and $\alpha_0,\alpha_\IF$ are defined in {\bf AI-AII}.}
\begin{itemize}
\item $\varphi=\lim_{u\to\infty}\frac{T_u}{u^{1/(2\alpha_\infty)}}$
\item $a_1=\frac{1+c\varphi}{\sqrt{2}A_\IF\varphi}$
\item $a_2=\frac{(1+c\varphi)^2}{A_\IF\varphi^2}\left(\alpha_\IF-\frac{c\varphi}{1+c\varphi}\right)$
\item $A=\left(\frac{{\alpha_\IF}}{c(1-{\alpha_\IF})}\right)^{-\alpha_\IF}\frac{1}{1-\alpha_\IF}$
\item $B=\left(\frac{\alpha_\IF}{c(1-\alpha_\IF)}\right)^{-\alpha_\IF-2}\alpha_\IF$
\item $m(u,t)=Var^{-1/2}\left(  \frac{X(t)}{u+ct} \right)= \frac{u+ct}{\sigma(t)}$
\item $\kd{t_u}=\arg \min_{t \ge 0} m(u,t)$
\item $\kd{t^*}=\frac{\alpha_\infty}{c(1-\alpha_\infty)}$
\item $\gamma=\lim_{u\to\infty}\frac{T_u}{u}$
\item $\Delta(u,s)=
\overleftarrow{\sigma}
 \left(\frac{\sqrt{2}\sigma^2(s)}{u+cs}\right)$
\item $\Omega(u,t)=\frac{m^2(u,t)\Delta(u,t)}{t}$
\item $\mu_\varphi(t)=\left\{\begin{array}{cc}
B_{\alpha_0}(t), & \hbox{if } \varphi=0\\
 \frac{1+c{\gamma}}{\sqrt{2}A_\IF\varphi^{2\alpha_\IF}}X
\left(\overleftarrow{\sigma}\left(\frac{\sqrt{2}A_\IF\varphi^{2\alpha_\IF}}{1+c{\gamma}}\right)t\right), & \hbox{if } \varphi\in (0,\IF)\\
B_{\alpha_\IF}(t), & \hbox{if } \varphi=\IF
\end{array}\right.$
\item $\eta_{\alpha_\IF}(t)=\left\{
            \begin{array}{ll}
B_{\alpha_0}(t), & \hbox{if } \alpha_\IF<1/2 \\
\frac{1+ct^*}{\sqrt{2}A_\IF t^*}X\left(\overleftarrow{\sigma}\left(\frac{\sqrt{2}A_\IF t^*}{1+ct^*}\right)t\right), & \hbox{if }\alpha_\IF=1/2\\
B_{\alpha_\IF}(t),&  \hbox{if } \alpha_\IF>1/2
              \end{array}
            \right.$
\end{itemize}

\bl{In the following lemma we give a
version of Theorem 3.5 in \cite{KEP2016}.}
\BEL\label{PPTH0}
%Let $X_u(t), t\in[a(u), b(u)]$ with $0 \in [a(u), b(u)]$,  be a family of centered continuous
%Gaussian processes satisfying (\ref{vvar2}) and (\ref{ccor2}).\\
Let  $X_u(t), t\in [a(u), b(u)]$ with $0\in [a(u), b(u)]$ be a family of centered continuous
Gaussian \rd{processes} with variance function $\sigma_u^2(t)$
satisfying, {as $u\rw\IF$},
\BQN\label{vvar2}
\sigma_u(0)=1, \ \ 1-\sigma_u(t)\sim \frac{|t|^{\beta}}{g(u)}, \quad  t\in [a(u), b(u)],
\EQN
with $\beta>0$, $\lim_{u\rw\IF}g(u)=\IF$, $\lim_{u\rw\IF}\frac{|a(u)|^{\beta}+|b(u)|^{\beta}}{g(u)}=0$,
and  correlation function satisfying
\BQN\label{ccor2}
\lim_{u\rw\IF}\sup_{s,t\in [a(u),b(u)], s\neq t}\left|n^2(u)\frac{1-Corr(X_u(s), X_u(t))}{\frac{Var(\eta(\Delta(u)|s-t|))}{Var(\eta(\Delta(u)))}}-1\right|=0,
\EQN
with
$\lim_{u\rw\IF}n(u)=\IF$ and
$\lim_{u\rw\IF}\Delta(u)=\theta\in [0,\IF],$
\bl{where  $\eta$ is a centered continuous Gaussian process with stationary increments, $\eta(0)=0$ and variance function satisfying {\bf AI-AII}.}
Let
\BQNY\label{vf} V_{\theta}(t)=\left\{\begin{array}{cc}
B_{\alpha_{0}}(t), & \theta=0\\
\frac{1}{\sqrt{Var(\eta(\theta))}}\eta(\theta t), & \theta\in (0,\IF) \\
B_{\alpha_{\IF}}(t), & \theta=\IF.
\end{array}\right.
\EQNY
{Suppose that} $\lim_{u\rw\IF}\frac{n^2(u)}{g(u)}=\nu\in [0,\IF]$.\\
i)If  $\nu=0$ and  $\lim_{u\rw\IF}\frac{(n(u))^{2/\beta}a(u)}{(g(u))^{1/\beta}}=y_1, \quad \lim_{u\rw\IF}\frac{(n(u))^{2/\beta}b(u)}{(g(u))^{1/\beta}}=y_2,
\quad  \lim_{u\rw\IF}\frac{(n(u))^{2/\beta}(a^2(u)+b^2(u))}{(g(u))^{2/\beta}}=0,$
 with $-\IF\leq y_1<y_2\leq \IF$,   then
\BQNY
\pk{\sup_{t\in [a(u), b(u)]}X_u(t)>n(u)}
\sim\mathcal{H}_{V_{\theta}}\int_{y_1}^{y_2}e^{-|s|^{\beta}}ds \left(\frac{g(u)}{n^2(u)}\right)^{1/\beta}\Psi(n(u)).
\EQNY
 ii) If $\nu\in (0,\IF)$ and $\lim_{u\rw\IF}a(u)=a\in[-\IF,0], \lim_{u\rw\IF}b(u)=b\in [0,\IF] $, then
 \BQNY
\pk{\sup_{t\in [a(u), b(u)]}X_u(t)>n(u)}
\sim
\mathcal{P}_{V_{\theta}}^{h}[a,b]\Psi(n(u)),
\EQNY
where $h(t)=\nu|t|^{\beta}$. \\
iii) If $\nu=\IF$,  then
\BQNY
\pk{\sup_{t\in [a(u), b(u)]}X_u(t)>n(u)}
 \sim \Psi(n(u)).
\EQNY
\EEL

We next focus on the analysis of  the behavior of variance and correlation
functions of the related Gaussian processes and Gaussian fields.
Hereafter, let $\overline{X}:=\frac{X}{\sqrt{Var(X)}}$ and   denote by $\dot{h}$
and $\ddot{h}$ the first and second derivative of \kd{twice continuously
differentiable function} $h$ respectively. Furthermore, for $X$ being a Gaussian processes with stationary increments satisfying {\bf AI-AII}, set \BQN\label{r_u}
r_{u}(s,t)=\E{\frac{X(ut)}{\rd{\sigma(ut)}} \frac{X(us)}{\sigma(us)}}, \quad s\neq t.\EQN
Suppose for a while that {\bf T1} holds. Then
\BQN\label{tr1}
\pi_{0,T_u}(u)&=&\pk{\sup_{t\in [0,T_u]}X(t)-ct>u}\nonumber\\
&=&\pk{\sup_{t\in [0,1]}X(T_ut)-c T_u t>u}\nonumber\\
&=&\pk{\sup_{t\in [0,1]}\frac{X(T_ut)}{u+cT_ut}>1}\nonumber\\
&=& \pk{\sup_{t\in [0,1]}\frac{X(T_ut)}{u+cT_ut}m(u,T_u)>m(u,T_u)}.
\EQN
Analogously, if {\bf T2} is satisfied, then
\BQNY\label{tr2}
\pi_{0,T_u}(u)= \pk{\sup_{t\in [0,T_u/u]}\frac{X(ut)}{u(1+ct)}m(u,t_u)>m(u,t_u)}.
\EQNY
Recall that $t^*=\frac{\alpha_\IF}{c(1-\alpha_\IF)}$ and $t_u=\arg\min_{t\geq 0} m(u,t)$.

\BEL\label{L1}
i)  Suppose that {\bf T1} and {\bf AI} are satisfied. Then for $u$ sufficiently large, the unique minimizer of $m(u, \cdot)$ over $[0, T_u]$ is $T_u$. Moreover, for each $u$ sufficiently large,
\BQNY\label{V4}
\frac{m(u,T_u)}{m(u, T_u t)}=1-a_u(1-t)(1+o(1)),\ \ \ t\rw 1,
\EQNY
where $a_u\rw \alpha_\IF-\frac{c\gamma}{1+c\gamma}$.\\
ii)
Suppose that {\bf AI} is satisfied.
 For $u$ large enough $t_u$ is unique, and  $t_u/u\rw t^*$, as
 $u\to \infty$, and $m(u,\cdot)$ is increasing over $[t_u,\IF)$. Moreover,  for each $u$ sufficiently large,
\BQNY\label{V4}
\frac{m(u,t_u)}{m(u,ut)}=1-b_u(t-t_u/u)^2(1+o(1)),\ \ \ t\rw t_u/u,
\EQNY
where $b_u\rw\frac{B}{2A}$ with $A, B$ defined in (\ref{AB}).
\EEL
\prooflem{L1}
Since the proofs of case i) and ii) are similar, we focus on detailed derivations only for case i)(see also Lemma 3.3 in \cite{KrzysPeng2015} for the proof of case ii)).\\
 We first note that for  $u$ sufficiently large, the minimizer of $m(u,\cdot)$ over $[0,T_u]$ is larger than any positive constant $T$.  Thus we focus on the the interval $[T,T_u]$. Theorem 1.7.2 in \cite{BI1989} yields that
\BQNY
\frac{(\dot{\sigma}(t))^2-\ddot{\sigma}(t)\sigma(t)}{(\dot{\sigma}(t))^2}-1&=&-\frac{\ddot{\sigma}(t)}{\dot{\sigma}(t)}
\frac{\sigma(t)}{\dot{\sigma}(t)}\\
&=&-\frac{2\sigma^2(t)}{t\dot{\sigma^2}(t)}\left(\frac{t\ddot{\sigma^2}(t)}{\dot{\sigma^2}(t)}
-\frac{t\dot{\sigma^2}(t)}{2\sigma^2(t)}\right)\rw \frac{1-\alpha_\IF}{\alpha_\IF}>0, \quad  t\rw\IF,
\EQNY
which implies that
$\frac{\sigma(t)}{\dot{\sigma}(t)}-t$ is increasing on interval $[T,\IF]$ for $T$ large enough. Further, Theorem 1.7.2 in \cite{BI1989} leads to, for $T$ sufficiently large and $u\rw\IF$,
\BQNY
\frac{\sigma(t)}{\dot{\sigma}(t)}-t-\frac{u}{c}\leq \frac{\sigma(T_u)}{\dot{\sigma}(T_u)}-T_u-\frac{u}{c}&=&\frac{u}{c}\left(\frac{\sigma(T_u)}{T_u\dot{\sigma}(T_u)}c\frac{T_u}{u}-c\frac{T_u}{u}-1\right)\\
&\sim& \frac{u}{c}\left(\frac{c(1-\alpha_\IF)}{\alpha_\IF}\gamma-1\right)<0, \quad  t\in [T,T_u],
\EQNY
and as $t\rw\IF$,
$$\dot{\sigma}(t)=\left(\frac{\sigma(t)}{t}\right)\left(\frac{t\dot{\sigma}(t)}{\sigma(t)}\right)\sim\frac{\alpha_\IF\sigma(t)}{t}>0.$$
Therefore,
\BQNY
\dot{m}(u, t)=\frac{c\sigma(t)-\dot{\sigma}(t)(u+ct)}{\sigma^2(t)}=\frac{c\dot{\sigma}(t)}{\sigma^2(t)}\left(\frac{\sigma(t)}{\dot{\sigma}(t)}-t-\frac{u}{c}\right)<0, \ \ t\in [T,T_u],
\EQNY
which implies that the minimum point is unique and equal to $T_u$.
Moreover, by \rdd{Theorem 1.7.2 and uniform convergence theorem }in \cite{BI1989}  we have
\BQNY
1-\frac{m(u,T_u)}{m(u, T_ut)}&=&1-\frac{\sigma(T_ut)}{\sigma(T_u)}\frac{u+cT_u}{u+cT_ut}\\
&\sim&\frac{T_u\dot{\sigma}(T_u\theta)}{\sigma(T_u)}(1-t)+1-\frac{1}{1-\rd{\frac{cT_u/u}{1+cT_u/u}}(1-t)}\\
&\sim&\left(\alpha_\IF-\frac{c\gamma}{1+c\gamma}\right)(1-t), \ \ t\rw 1,
\EQNY
with $\theta\in (t,1)$.
This completes the proof. \QED
%The proof of Lemma \ref{L1} is postponed to Appendix.

\bl{In the following lemma}  we derive asymptotic behaviour of $r_u(s,t)$ and $r_{T_u}(s,t)$, needed while applying Lemma \ref{PPTH0}.
\BEL\label{L2}
Suppose that {\bf AI, AII} hold. Then for any  $T_u\rw\IF, \delta_u\rw 0$, as $u\rw \IF$,
 \BQNY\label{cor0}
\lim_{u\rw \IF}\sup_{s\neq t, |t-t_u/u|, |s-t_u/u|<\delta_u}\left|\frac{1-r_u(s,t)}{\frac{\sigma^2(u|s-t|)}{2\sigma^2(ut^*)}}-1\right|=0,
\EQNY
and
\BQNY\label{cor1}
\lim_{u\rw \IF}\sup_{s\neq t, s,t\in[1-\delta_u,1]}\left|\frac{1-r_{T_u}(s,t)}{\frac{\sigma^2(T_u|s-t|)}{2\sigma^2(T_u)}}-1\right|=0.
\EQNY
\EEL
\bl{ In the next lemma we collect some asymptotics which will be helpful in the proofs.}
\BEL\label{Formula}
Let $\mathbb{Q}_i>0, i=1,...,5$ be some  constants.\\
i) \bl{ If  {\bf T1} is satisfied then, as $u\rw\IF$,
$$\Omega(u,T_u)\sim \left\{\begin{array}{cc}
\mathbb{Q}_1\frac{u}{T_u}\left(\frac{T_u}{u^{1/(2\alpha_\IF)}}\right)^{\frac{2(1-\alpha_0)\alpha_\IF}{\alpha_0}}, & \varphi=0\\
\mathbb{Q}_2u^{1-1/(2\alpha_\IF)},& \varphi\in (0,\IF)\\
\mathbb{Q}_3\frac{u}{T_u}\left(\frac{T_u}{u^{1/(2\alpha_\IF)}}\right)^{2(1-\alpha_\IF)},& \varphi=\IF.
\end{array}\right.$$}
%\rd{where $\mathbb{Q}_i>0, i=1,2,3$ are some  constants.}\\
ii) \bl{ $\frac{m(u,t_u)\Delta(u,t_u)}{u}\sim \mathbb{Q}_4\frac{\overleftarrow{\sigma}(u^{-1}\sigma^2(u))}{\sigma(u)}\sim \mathbb{Q}_5u^\beta$, as $u\rw\IF$, with $\beta<0$ defined by
\BQN\label{beta}\beta=\left\{\begin{array}{cc}
\frac{2\alpha_\IF-1}{\alpha_0}-\alpha_\IF, & \text{if} \ \ \alpha_\IF<1/2\\
-\alpha_\IF, & \text{if} \ \ \alpha_\IF=1/2\\
\frac{2\alpha_\IF-1}{\alpha_\IF}-\alpha_\IF, & \text{if} \ \ \alpha_\IF>1/2.
\end{array}\right.\EQN}
\EEL
%where \rd{$\mathbb{Q}_i>0, i=4,5$ are two  constants and  $\beta<0$ for all $0<\alpha_\IF<1$}.\\
We conclude this section with the study of the  limit of $\Omega(u,T_u)$, which determines the asymptotics for the short-time horizon case.
\BEL\label{LA}
Assume  that {\bf T1} is satisfied and $\varphi\in [0,\IF]$. Then $\lim_{u\rw\IF}\Omega(u,T_u)$ exists and
\\ i)
if $\varphi=0$, then $\lim_{u\rw\IF}\Omega(u,T_u)
\in [0,\IF]$;\\
ii) if $\varphi\in (0,\IF)$, then $\lim_{u\rw\IF}\Omega(u,T_u)\in (0,\IF]$;\\
iii) if $\varphi=\IF$, then $\lim_{u\rw\IF}\Omega(u,T_u)=\IF$.\\
\EEL
The proofs of  Lemmas \ref{L2}-\ref{LA} are standard but need some tedious calculations; thus we skip the proofs
referring to related derivations in, e.g.,  \cite{KrzysPeng2015}.

%The proof of Lemma \ref{LA} is postponed to Appendix.

\section{Proofs of main results}
{In the rest of the paper,  by $\mathbb{Q}, \mathbb{Q}_i, i=1,2,3, \ldots$
we denote some positive constants that may differ from line to line.}
\label{s.proofs}
If multiple limits appear, we shall write \rdd{
$ f_u(S,S_1,\epsilon) \sim f^*(u),  u\to \IF, S\to \IF, S_1 \to \IF, \epsilon\rw 0$
to mean that
$$\lim_{\epsilon\rw 0}\lim_{S_1\to \IF} \lim_{S\rw\IF}\lim_{u\to \IF} \frac{f_u(S,S_1,\epsilon)}{f^*(u)}=1.$$}
\rdd{Let $f_\lambda(t)=\frac{\sigma^2(t)}{t^\lambda}, t>0$ with $\lambda\in (0, \min(2\alpha_0, 2\alpha_\IF))$. Following {\bf AI-AII}, $f_\lambda$ is a regularly varying function at $0$ and $\IF$ with index $2\alpha_0-\lambda$ and $2\alpha_\IF-\lambda$ respectively. For any $T>0$,  by uniform convergence theorem, e.g. \cite{BI1989}, we have that
$$\lim_{u\rw\IF}\sup_{t\in (0,T]}\left|\frac{f_\lambda(ut)}{f_\lambda(u)}-t^{2\alpha_\IF-\lambda}\right|=0,$$
implying that for $\lambda\in (0, \min(2\alpha_0, 2\alpha_\IF))$ and $u$ sufficiently large,
\BQN\label{reg}
\frac{\sigma^2(ut)}{\sigma^2(u)}=\frac{f_\lambda(ut)}{f_\lambda(u)}t^\lambda\leq 2T^{2\alpha_\IF-\lambda}t^\lambda, \quad t\in (0,T].
\EQN}

\subsection{Proof of Proposition \ref{p.T}}
Observe that
$$\pi_{0,T}(u)=\pk{\sup_{0\leq t\leq T}\frac{X(t)}{u+ct}\frac{u+cT}{\sigma(T)}>\frac{u+cT}{\sigma(T)}}.$$
One can easily check that
%assumption that
%$\dot{\sigma^2}(t)>0$ for $t\in (0,T]$,
$\sup_{0\leq t\leq T}Var\left(\frac{X(t)}{u+ct}\frac{u+cT}{\sigma(T)}\right)=1$.
Moreover, by {\bf AII}, there exists $C>0$ such that for all $0\leq s\leq t\leq T$ and \rd{$u>1$
\BQNY
\frac{(u+cT)^2}{\sigma^2(T)}\E{\left(\frac{X(t)}{u+ct}-\frac{X(s)}{u+cs}\right)^2}&=&\frac{(u+cT)^2}{\sigma^2(T)}
\E{\left(\frac{X(t)-X(s)}{u+ct}+\frac{c(s-t)X(s)}{(u+ct)(u+cs)}\right)^2}\\
&\leq& 2\frac{(u+cT)^2}{\sigma^2(T)}\left(\frac{\sigma^2(|t-s|)}{(u+ct)^2}+\frac{c^2\sigma^2(s)(t-s)^2}{(u+ct)^2(u+cs)^2}\right)\leq C|t-s|^{\alpha_0}.
\EQNY}
Hence by Piterbarg inequality ( Theorem 8.1 in \cite{Pit96}), we have for $u$ sufficiently large,
\begin{eqnarray}
\pi_{0,T}(u)\leq \mathbb{Q}\left(\frac{u+cT}{\sigma(T)}\right)^{2/\alpha_0}\Psi\left(\frac{u+cT}{\sigma(T)}\right)\label{Pit1},
\end{eqnarray}
 implying that, as $u\rw\IF$,
$$\pi_{0,T}(u)=o\left(\Psi\left(\frac{u-x+cT}{\sigma(T)}\right)\right).$$
Thus in view of (\ref{Basic1})
$$\pi_{x,T}(u)\sim \Psi\left(\frac{u-x+cT}{\sigma(T)}\right), \quad u\rw\IF.$$
This completes the proof.\QED

\subsection{Proof of Theorem \ref{submain1}}
The idea of the proof is based on the observation, by (\ref{tr1}), that
\BQN\label{e1}
\Pi_1(u)\leq \pi_{0, T_u}(u)\leq \Pi_1(u)+ \Pi_2(u),
\EQN
where \BQNY
\Pi_1(u)&=&\pk{\sup_{t\in E(u)}\frac{X(T_ut)}{u+cT_ut}m(u,T_u)>m(u,T_u)}, \quad E(u)=[1-\left((\ln m(u,T_u))/m(u,T_u)\right)^2, 1],\\
\Pi_2(u)&=&\pk{\sup_{t\in [0,1]\setminus E(u)}\frac{X(T_ut)}{u+cT_ut}m(u,T_u)>m(u,T_u)}.
\EQNY
\gr{In what follows, we shall derive the exact asymptotics of $\Pi_1(u)$ by applying  Lemma \ref{PPTH0} and then show that $\Pi_2(u)$ is asymptotically negligible compared with $\Pi_1(u)$ as $u\rw\IF$.}\\
{\it\underline{Analysis of $\Pi_1(u)$}}. In order to apply Lemma \ref{PPTH0} , we rewrite
$$\Pi_1(u)=\pk{\sup_{t\in [0, b(u)]}X_u(t)>m(u,T_u)}, $$
where
 $$X_u(t):=\frac{X(T_u-\Delta(u,T_u)t)}{u+cT_u-c\Delta(u,T_u)t}m(u,T_u), t\in [0, b(u)], \quad \hbox{with } b(u)=\left((\ln m(u,T_u))/m(u,T_u)\right)^2\frac{T_u}{\Delta(u,T_u)}.$$  Let $\sigma_u(t)=\sqrt{Var(X_u(t))}, \quad  g(u)=\frac{T_u}{a_u\Delta(u,T_u)},$
 with $a_u$ defined in Lemma \ref{L1}.

In light of Lemma \ref{L1} combined with Lemma \ref{L2}, we have
$$
\sigma_u(0)=1, \ \ 1- \sigma_u(t)\sim \frac{t}{g(u)}, t\in [0, b(u)], \ \ \lim_{u\rw \IF}\sup_{s,t \in[0, b(u)], s\neq t}\left|\frac{m^2(u,T_u)(1-Corr(X_u(s), X_u(t)))}{\frac{\sigma^2(\Delta(u,T_u)|s-t|)}{\sigma^2(\Delta(u,T_u))}}-1\right|=0.
$$
Moreover,
\BQNY\label{neq4}
\nu=\lim_{u\rw\IF} \frac{m^2(u,T_u)}{g(u)}=\left(\alpha_\IF-\frac{c\gamma}{1+c\gamma}\right)\lim_{u\rw\IF}\Omega(u,T_u),\ \ \lim_{u\rw\IF}g(u)=\lim_{u\rw\IF}m(u,T_u)=\IF, \quad \lim_{u\rw\IF}\frac{b(u)}{g(u)}=0.
\EQNY
Hence, if  $\lim_{u\rw\IF}\Omega(u,T_u)=0$, then using that
\BQNY\label{neq5}
 y_1=0,\quad  y_2=\lim_{u\rw\IF}\frac{b(u)m^2(u,T_u)}{g(u)}=\IF, \ \ \lim_{u\rw\IF}\left(\frac{b(u)}{g(u)}\right)^2m^2(u,T_u)=0,
\EQNY
i) of Lemma \ref{PPTH0} leads to, as $u\rw\IF$,
\BQNY
\Pi_1(u)&\sim& \mathcal{H}_{B_{\alpha_0}}\int_{y_1}^{y_2} e^{-t}dt\left(\alpha_\IF-\frac{c\gamma}{1+c\gamma}\right)^{-1}
\Omega^{-1}(u,T_u)
\Psi(m(u,T_u))\\
&\sim& \mathcal{H}_{B_{\alpha_0}}\int_0^\IF e^{-t}dt\left(\alpha_\IF-\frac{c\gamma}{1+c\gamma}\right)^{-1}
\Omega^{-1}(u,T_u)
\Psi(m(u,T_u))\\
&\sim& \mathcal{H}_{B_{\alpha_0}}\left(\alpha_\IF-\frac{c\gamma}{1+c\gamma}\right)^{-1}
%\frac{T_u}{\Delta(u,T_u)m^2(u,T_u)}
\Omega^{-1}(u,T_u)
\Psi(m(u,T_u)).
\EQNY}
If $\lim_{u\rw\IF}\Omega(u,T_u)=\lim_{u\rw\IF}\Omega_\IF\in (0,\IF)$, then $\nu=\left(\alpha_\IF-\frac{c\gamma}{1+c\gamma}\right)\Omega_\IF$ and \BQNY\label{nu}\lim_{u\rw\IF}\Delta(u,T_u)=
\rd{\lim_{u\rw\IF}}\overleftarrow{\sigma}
 \left(\frac{\sqrt{2}\sigma^2(T_u)}{u+cT_u}\right)=\left\{\begin{array}{cc}
 0, &\hbox{if } \varphi=0\\
\overleftarrow{\sigma}(\frac{\sqrt{2}A_\IF\varphi^{2\alpha_\IF}}{1+c{\gamma}}), &\hbox{if } \varphi\in (0,\IF).
 \end{array}\right.\EQNY
 Note that due to Lemma \ref{LA}, we can exclude  case $\varphi=\IF$. Thus, \gr{by case ii}) in Lemma \ref{PPTH0}, we have
\BQNY
\Pi_1(u)\sim \mathcal{P}_{\mu_\varphi}^{\Omega_\IF(\alpha_\IF-\frac{c\gamma}{1+c\gamma})t}
\Psi(m(u,T_u)).
\EQNY
If $\lim_{u\rw\IF}\Omega(u,T_u)=\IF$, \gr{in light of case iii) }in Lemma \ref{PPTH0}, we have
\BQNY
\Pi_1(u)\sim
\Psi(m(u,T_u)).
\EQNY
{\it \underline{Analysis of $\Pi_2(u)$}}.
Due to Lemma \ref{L1}, we have
$$\sup_{t\in [0,1]\setminus E(u)} Var\left(\frac{X(T_ut)}{u+cT_ut}m(u,T_u)\right)\leq 1-\mathbb{Q}(m(u,T_u))^{-2}(\ln m(u,T_u))^{2}. $$
Moreover, \rdd{by (\ref{reg})} we have
$$\E{\left(\frac{X(T_ut)}{u+cT_ut}m(u,T_u)-\frac{X(T_us)}{u+cT_us}m(u,T_u)\right)^2}\leq \mathbb{Q}\left(\frac{\sigma^2(T_u|t-s|)}{\sigma^2(T_u)}+|t-s|^2\right)\leq \mathbb{Q}|t-s|^\lambda,$$
with $0<\lambda<\min (2\alpha_0, 2\alpha_\IF).$
Thus, by %Lemma 5.1 in \cite{KEP2015}
Piterbarg inequality (Theorem 8.1 in \cite{Pit96}), we have
\BQNY\label{ad1}
\Pi_2(u)&\leq&\mathbb{Q} (m(u,T_u))^{2/\lambda}\Psi\left(\frac{m(u,T_u)}{\sqrt{1-\mathbb{Q}(m(u,T_u))^{-2}(\ln m(u,T_u))^{2}}}\right)\\
&\leq& \mathbb{Q} (m(u,T_u))^{2/\lambda}e^{-\mathbb{Q} (\ln m(u,T_u))^2}\Psi\left( m(u,T_u))\right)\\
&=&o\left(\Pi_1(u)\right), \quad u\rw\IF.
\EQNY
Hence, $\pi_{0,T_u}(u)\sim \Pi_1(u)$, as $u\rw\IF$, which
completes the proof.\QED

\subsection{Proof of Theorem \ref{nth1}}
Recall that
 \BQN\label{nneq10}
\pi_{x,T_u}(u)&=&\pi_{0,T_u}(u)+\Psi\left(\frac{u-x+cT_u}{\sigma(T_u)}\right)\nonumber\\
&& \ \ -\pl{X(T_u)-c T_u>u-x, \sup_{0\leq s\leq T_u}\left(X(T_u)-X(s)-c(T_u-s)\right)>u}.
\EQN
 Hence, if $\pi_{0,T_u}(u)$ or $\Psi\left(\frac{u-x+cT_u}{\sigma(T_u)}\right)$ is asymptotically dominating, then it determines the  asymptotics of $\pi_{x,T_u}(u)$ as $u\rw\IF$. By Theorem \ref{submain1}, we have
\BQN\label{ratio}
\frac{\pi_{0,T_u}(u)}{\Psi\left(\frac{u-x+cT_u}{\sigma(T_u)}\right)}\sim\left\{\begin{array}{cc}
\mathcal{H}_{B_{\alpha_0}}\left(\alpha_\IF-\frac{c\gamma}{1+c\gamma}\right)^{-1}
\Omega^{-1}(u,T_u)e^{-\frac{u+cT_u}{\sigma^2(T_u)}x},& \text{if} \ \ \lim_{u\rw\IF}\Omega(u,T_u)=0 \\
\mathcal{P}_{\mu_\varphi}^{\Omega_\IF(\alpha_\IF-\frac{c\gamma}{1+c\gamma})t}e^{-\frac{u+cT_u}{\sigma^2(T_u)}x}, & \text{if} \ \ \lim_{u\rw\IF}\Omega(u,T_u)=\Omega_\IF\in (0,\IF) \\
e^{-\frac{u+cT_u}{\sigma^2(T_u)}x},& \text{if} \ \ \lim_{u\rw\IF}\Omega(u,T_u)=\IF. \end{array}\right.
\EQN
We shall distinguish three cases related to the value of $\varphi$.\\
{$\diamond$ \underline{Case $\varphi=0$}}. \bl{Observe that,  for $u$ sufficiently large,
\BQNY\label{bas}
\frac{u+cT_u}{\sigma^2(T_u)}\geq \mathbb{Q}\left(\frac{T_u}{u^{1/(2\alpha_\IF)}}\right)^{-2\alpha_\IF}.
\EQNY}
Moreover, \rd{for  $u$ sufficiently large,} i) \bl{of Lemma \ref{Formula} leads to
$$\Omega^{-1}(u,T_u)\leq \mathbb{Q}_1\frac{T_u}{u}\left(\frac{T_u}{u^{1/(2\alpha_\IF)}}\right)^{-\frac{2(1-\alpha_0)\alpha_\IF}{\alpha_0}}\leq \mathbb{Q}_1\left(\frac{T_u}{u^{1/(2\alpha_\IF)}}\right)^{-\frac{2(1-\alpha_0)\alpha_\IF}{\alpha_0}}.$$
The above implies that
$$e^{-\frac{u+cT_u}{\sigma^2(T_u)}x}\leq e^{-\mathbb{Q}\left(\frac{T_u}{u^{1/(2\alpha_\IF)}}\right)^{-2\alpha_\IF}}\rw 0, \quad u\rw\IF,$$
and
$$\Omega^{-1}(u,T_u)e^{-\frac{u+cT_u}{\sigma^2(T_u)}x}\leq \mathbb{Q}_1\left(\frac{T_u}{u^{1/(2\alpha_\IF)}}\right)^{-\frac{2(1-\alpha_0)\alpha_\IF}{\alpha_0}}e^{-\mathbb{Q}\left(\frac{T_u}{u^{1/(2\alpha_\IF)}}\right)^{-2\alpha_\IF}}\rw 0, \quad u\rw\IF.$$
Hence, by (\ref{ratio})}, we have that
$$\pi_{0,T_u}(u)=o\left(\Psi\left(\frac{u-x+cT_u}{\sigma(T_u)}\right)\right), \ \ u\rw\IF,$$
which implies that
$\pi_{x,T_u}\sim \Psi\left(\frac{u-x+cT_u}{\sigma(T_u)}\right).$\\
{$\diamond$ \underline{Case $\varphi\in (0,\IF)$}}. Note  that
\bl{$$\lim_{u\rw\IF}\frac{u+cT_u}{\sigma^2(T_u)}=\frac{1+c\gamma}{A_\IF\varphi^{2\alpha_\IF}}\in (0,\IF)$$}
and $$\lim_{u\rw\IF}\Omega(u,T_u)\in (0,\IF],$$
 by Lemma \ref{LA}. Thus  (\ref{ratio}) leads to $$\pi_{0,T_u}(u)=O\left(\Psi\left(\frac{u-x+cT_u}{\sigma(T_u)}\right)\right), \quad u\rw\IF.$$
Hence, this case needs another approach than applied in scenario $\varphi=0$.
We have  $T_u\sim \varphi u^{1/(2\alpha_\IF)}$ and $\lim_{u\rw\IF}\frac{T_u}{u}<t^*$,  \rd{so}  $\alpha_\IF\geq 1/2$. \bl{Moreover,}
i) \bl{ of Lemma \ref{Formula} gives that
$$\Omega(u,T_u)\sim \mathbb{Q}_1 u^{1-1/(2\alpha_\IF)}, \quad u\rw\IF.$$}
{$\diamond$ \underline{Subcase $\varphi\in (0,\IF)$, $\alpha_\IF=1/2$}}.
If $\alpha_\IF=1/2$, then $\lim_{u\rw\IF}\Omega(u,T_u)=\Omega_\IF\in (0,\IF)$.
Recall that
$$\pi_{x,T_u}=\pl{X(T_u)-c T_u>u-x \hbox{ or } \sup_{0\leq s\leq T_u}\left(X(T_u)-X(s)-c(T_u-s)\right)>u},$$
and observe that for any $S>0$,
\BQN\label{nin1}
\Pi_x^{(1)}(u)\leq \pi_{x,T_u}(u)\leq \Pi_x^{(1)}(u)+\Pi^{(2)}(u),
\EQN
where
\BQNY
\Pi_x^{(1)}(u)&=&\pl{X(T_u)-c T_u>u-x \hbox{ or } \sup_{0\leq s\leq S}\left(X(T_u)-X(s)-c(T_u-s)\right)>u}, \\
 \Pi^{(2)}(u)&=&\pk{\sup_{S\leq s\leq T_u}\left(X(T_u)-X(s)-c(T_u-s)\right)>u}.
\EQNY
\gr{Next we shall first derive the exact asymptotics of $\Pi_x^{(1)}(u)$ and then show that $\Pi^{(2)}(u)=o\left(\Pi_x^{(1)}(u)\right)$ as $u\rw\IF, S\rw\IF$.}\\
{\it\underline{Analysis of $\Pi_x^{(1)}(u)$}}.
The same transformation as given in (\ref{tr1}) leads to
$$\Pi_x^{(1)}(u)=\pk{X(T_u)>u-x+cT_u \ \ \text{or} \ \ \sup_{t\in I_0(u)}\left(\frac{X(T_u)-X(T_ut)}{u+cT_u(1-t)}\right)m(u,T_u)>m(u,T_u)}, \quad I_0(u)=\left[0,\frac{S}{T_u}\right].$$
 Using Lemma \ref{L1}, we have for $0<\epsilon<\alpha_\IF-\frac{c\gamma}{1+c\gamma}$ and $u$ large enough,
 \BQNY
\Pi_x^{1,+\epsilon}(u)\leq \Pi_x^{(1)}(u)\leq \Pi_x^{1,-\epsilon}(u)
 \EQNY
with
$$\Pi_x^{1,\pm\epsilon}(u)=\pk{\overline{X}(T_u)>\frac{u-x+cT_u}{\sigma(T_u)} \ \ \text{or} \ \ \sup_{t\in I_{0}(u)}\frac{\overline{X(T_u)-X(T_ut)}}{1+(\alpha_\IF-\frac{c\gamma}{1+c\gamma}\pm\epsilon)t}>m(u,T_u)}.$$
We first focus on $\Pi_x^{1,-\epsilon}(u)$. Let
$$g_u(w):=\pk{\overline{X}(T_u)>\frac{u-x+cT_u}{\sigma(T_u)} \ \ \text{or} \ \ \sup_{t\in I_{0}(u)}\frac{\overline{X(T_u)-X(T_ut)}}{1+(\alpha_\IF-\frac{c\gamma}{1+c\gamma}-\epsilon)t}>m(u,T_u)\Bigl\lvert\overline{X}(T_u)=m(u,T_u)-\frac{w}{m(u,T_u)}}.$$
Then
\BQN\label{newadd}
\Pi_x^{1,-\epsilon}(u)=\frac{1}{\sqrt{2\pi}m(u,T_u)}e^{-\frac{m^2(u,T_u)}{2}}\int_{-\IF}^{\IF}e^{w-\frac{w}{m^2(u,T_u)}}g_u(w)dw.
\EQN
Notice that if $w<\sqrt{2}a_1x-\rdd{\nu}$ with $\rdd{\nu>0}$ and  $a_1=\frac{1+c\varphi}{\sqrt{2}A_\IF\varphi}$, then for $u$ sufficiently large,
$$m(u,T_u)-\frac{w}{m(u,T_u)}>\frac{u-x+cT_u}{\sigma(T_u)},$$
implying that $g_u(w)=1$. Analogously,   if $w>\sqrt{2}a_1x+\rdd{\nu}$ with $\rdd{\nu>0}$, for $u$ sufficiently large, then
$$m(u,T_u)-\frac{w}{m(u,T_u)}<\frac{u-x+cT_u}{\sigma(T_u)},$$
which means that
$g_u(w)=\pk{\sup_{t\in I_{0}(u)}\frac{\overline{X(T_u)-X(T_ut)}}{1+(\alpha_\IF-\frac{c\gamma}{1+c\gamma}-\epsilon)t}>m(u,T_u)
\Bigl\lvert\overline{X}(T_u)=m(u,T_u)-\frac{w}{m(u,T_u)}}.$\\
\rd{In order to analyze the conditional process, let
\BQNY
R_{u}(s,t)=\E{\frac{X(T_u)-X(s)}{\sigma(T_u-s)}\frac{X(T_u)-X(t)}{\sigma(T_u-t)}}, \quad s,t\in [0,S].
\EQNY
Then by Taylor formula we have
\BQNY
1-R_u(s,t)&=&\frac{\sigma^2(|t-s|)-(\sigma(T_u-t)-\sigma(T_u-s))^2}{2\sigma(T_u-t)\sigma(T_u-s)}\nonumber\\
&=&\frac{\sigma^2(|t-s|)-(\dot{\sigma}(T_u-\theta))^2 (t-s)^2}{2\sigma(T_u-t)\sigma(T_u-s)},
\EQNY
where $\theta\in (s,t)$.
Moreover, Theorem 1.7.2 in \cite{BI1989} yields that
$$\dot{\sigma}(T_u-\theta)\sim \alpha_\IF \frac{\sigma(T_u)}{T_u}\sim \alpha_\IF \sqrt{A_\IF}T_u^{\alpha_\IF-1}, \quad u\rw\IF,$$
 and {\bf AII} leads to
$$\sigma^2(|t-s|)\geq \mathbb{Q}|t-s|^{2\alpha_0}, \quad s\neq t, s,t\in [0,S].$$
Hence
\BQNY
\sup_{s\neq t, s,t\in [0,S]}\frac{(\dot{\sigma}(T_u-\theta))^2 (t-s)^2}{\sigma^2(|t-s|)}\leq \mathbb{Q}\sup_{s\neq t, s,t\in [0,S]} T_u^{2(\alpha_\IF-1)}|t-s|^{2(1-\alpha_0)}\rw 0, \quad u\rw\IF.
\EQNY
The above analysis implies that
\BQN\label{Ru}
1-R_u(s,t)
\sim \frac{\sigma^2(|t-s|)}{2\sigma^2(T_u)},  s,t\in [0,S], \quad u\rw\IF.
\EQN
Notice that
$$\overline{X(T_u)-X(t)}=\left(\overline{X(T_u)-X(t)}-R_u(0,t)\overline{X}(T_u)\right)+R_u(0,t)\overline{X}(T_u),$$
where  $\overline{X(T_u)-X(t)}-R_u(0,t)\overline{X}(T_u), t\in [0,S]$ is independent of $\overline{X}(T_u)$.}
Thus
\BQNY
&&\left\{m(u,T_u)\frac{\overline{X(T_u)-X(t)}}{1+(\alpha_\IF-\frac{c\gamma}{1+c\gamma}-\epsilon)t/T_u}-m^2(u,T_u)+
w,  t\in [0,S]\Bigl\lvert\overline{X}(T_u)=m(u,T_u)-\frac{w}{m(u,T_u)}\right\}\\
 && \quad \stackrel{d}{=} \left\{\frac{Z_u(t)+h_u(w,t)}{1+(\alpha_\IF-\frac{c\gamma}{1+c\gamma}-\epsilon)t/T_u}, t\in [0,S]\right\},
\EQNY
where
$$Z_u(t)=m(u,T_u)\left(\overline{X(T_u)-X(t)}-R_u(0,t)\overline{X}(T_u)\right), \ \ t\in [0,S],$$
and
\BQN\label{hu}
h_u(w,t)&=&-m^2(u,T_u)(1-R_u(0,t))-\frac{m^2(u,T_u)}{T_u}\left(\alpha_\IF-\frac{c\gamma}{1+c\gamma}-\epsilon\right)t\nonumber\\
&& +w\left(1-R_u(0,t)
+\frac{\alpha_\IF-\frac{c\gamma}{1+c\gamma}-\epsilon}{T_u}
t\right).
\EQN
\rd{It follows that
\BQNY
Cov(Z_u(t), Z_u(s))&=&m^2(u,T_u)(R_u(s,t)-R_u(0,s)R_u(0,t))\nonumber\\
&=&m^2(u,T_u)\left((1-R_u(0,s))+(1-R_u(0,t))-(1-R_u(s,t))\right)\nonumber\\
&&\quad -m^2(u,T_u)(1-R_u(0,s))(1-R_u(0,t)), \quad  s,t\in [0,S].
\EQNY
Using (\ref{Ru}), we have that
\BQNY
\lim_{u\rw\IF}m^2(u,T_u)(1-R_u(s,t))=\lim_{u\rw\IF}\frac{m^2(u,T_u)}{2\sigma^2(T_u)}\sigma^2(|t-s|)=a_1^2\sigma^2(|t-s|),
\EQNY
where $a_1$ is defined in (\ref{a12}). Consequently,
\BQN\label{limvar}
\lim_{u\rw\IF}Cov(Z_u(t), Z_u(s))=Cov(\sqrt{2}a_1X(t), \sqrt{2}a_1X(s)), \quad s,t\in [0,S],
\EQN
\rdd{and for each $w\in \mathbb{R}$,}}
\BQNY
h_u(w,t)&\rw& -a_1^2\sigma^2(t)-a_2(\epsilon)t,\label{trend}
\EQNY
uniformly with respect to $t\in [0,S]$ with  $a_2(\epsilon)=\frac{(1+c\varphi)^2}{A_\IF\varphi^2}\left(\alpha_\IF-\frac{c\varphi}{1+c\varphi}-\epsilon\right)$. Moreover, for $u$ sufficiently large
\BQN\label{xeq10}
\E{\left(Z_u(t)-Z_u(s)\right)^2}
 &\leq& m^2(u,T_u)\left(\E{\left(\overline{X(T_u)-X(t)}-\overline{X(T_u)-X(s)}\right)^2}+(R_u(0,t)-R_u(0,s))^2\right)\nonumber\\
& \leq& 4m^2(u,T_u)(1-R_u(s,t))\nonumber\\
&\leq&  \mathbb{Q}\sigma^2(|t-s|) \leq \mathbb{Q}|t-s|^{\alpha_0/2}, \ \ s,t\in [0,S].
\EQN
Thus  $\{Z_u(t)+h_u(w,t), t\in [0,S]\}$ weakly converges to $\{\sqrt{2}a_1X(t)-a_1^2\sigma^2(t)-a_2(\epsilon)t, t\in[0,S]\}$.
Since $1+(\alpha_\IF-\frac{c\gamma}{1+c\gamma}-\epsilon)t/T_u$ uniformly converges to $1$ with respect to $t\in [0,S]$, then
$\left\{\frac{Z_u(t)+h_u(w,t)}{1+(\alpha_\IF-\frac{c\gamma}{1+c\gamma}-\epsilon)t/T_u}, t\in [0,S]\right\}$ weakly converges to $\left\{\sqrt{2}a_1X(t)-a_1^2\sigma^2(t)-a_2(\epsilon)t, t\in[0,S]\right\}$, implying that, as $u\rw\IF$,
$$g_u(w)\rw \pk{\sup_{t\in [0,S]}\left(\sqrt{2}a_1X(t)-a_1^2\sigma^2(t)-a_2(\epsilon)t\right)>w}:=g_\IF(w), \ \ w>\sqrt{2}a_1x+\rdd{\nu}.$$
 \rd{Noting that for  $u$  large enough and $w>0$
 $$\sup_{t\in [0,S]}h_u(w,t)\leq w\sup_{t\in [0,S]}\left(1-R_u(0,t)
+\frac{\alpha_\IF-\frac{c\gamma}{1+c\gamma}-\epsilon}{T_u}
t\right)\leq w/4,$$
 it follows that for $u$ large enough and $w>0$
\BQNY
g_u(w)&=&\pk{\sup_{t\in [0,S]}\frac{Z_u(t)+h_u(w,t)}{1+(\alpha_\IF-\frac{c\gamma}{1+c\gamma}-\epsilon)t/T_u}>w}\\
&\leq&\pk{\sup_{t\in [0,S]}Z_u(t)>w/2-\sup_{t\in [0,S]}h_u(w,t)}\\
&\leq&\pk{\sup_{t\in [0,S]}Z_u(t)>w/4}.
\EQNY}
Thus, in view of (\ref{limvar}) and (\ref{xeq10}), by Piterbarg inequality (Theorem 8.1 in \cite{Pit96}) we have for $w$ and $u$ sufficiently large,
\BQNY
\pk{\sup_{t\in [0,S]}Z_u(t)>w/4}\leq \mathbb{Q} w^{4/\alpha_0}\Psi\left(\frac{w}{8a_1\sigma(S)}\right).
\EQNY
\bl{Thus by  Mills' ratio we have that, for $u$ sufficiently large and $w>\sqrt{2}a_1x+\rdd{\nu}$,
\BQN\label{bound}
g_u(w)\leq \mathbb{Q} w^{4/\alpha_0}\Psi\left(\frac{w}{8a_1\sigma(S)}\right)
\leq \mathbb{Q}_1 w^{4/\alpha_0-1} e^{-\frac{w^2}{128a_1^2\sigma^2(S)}}.
\EQN}
Therefore the dominated convergence theorem leads to, as $u\rw\IF$,
$$\int_{\sqrt{2}a_1x+\rdd{\nu}}^\IF e^{w-\frac{w}{m^2(u,T_u)}}g_u(w)dw\rw \int_{\sqrt{2}a_1x+\rdd{\nu}}^\IF e^{w}g_\IF(w)dw.$$
Moreover,
$$\int_{-\IF}^{\sqrt{2}a_1x-\rdd{\nu}} e^{w-\frac{w}{m^2(u,T_u)}}g_u(w)dw\rw \int_{-\IF}^{\sqrt{2}a_1x-\rdd{\nu}} e^wdw, \ \ u\rw\IF,$$
and for $u$ sufficiently large,
$$\int_{\sqrt{2}a_1x-\rdd{\nu}}^{\sqrt{2}a_1x+\rdd{\nu}} e^{w-\frac{w}{m^2(u,T_u)}}g_u(w)dw\leq 2\rdd{\nu} e^{\sqrt{2}a_1x+\rdd{\nu}}.$$
Hence, \rd{due to (\ref{newadd}),}
$$\lim_{\rdd{\nu}\rw 0 }\limsup_{u\rw\IF}\frac{\Pi_x^{1,-\epsilon}(u)}{\Psi(m(u,T_u))}=\int_{-\IF}^{\sqrt{2}a_1x} e^wdw+ \int_{\sqrt{2}a_1x}^\IF e^{w}g_\IF(w)dw=\mathcal{P}_{a_1X}^{a_2(\epsilon)t, \sqrt{2}a_1x}[0,S].$$
 Similarly,
$$\liminf_{ u\rw\IF}\frac{\Pi_x^{1, +\epsilon}(u)}{\Psi(m(u,T_u))}= \mathcal{P}_{a_1X}^{a_{2}(-\epsilon)t, \sqrt{2}a_1x}[0,S]. $$
Thus, letting $\epsilon\rw 0$,
\BQN\label{new1}
\lim_{u\rw\IF}\frac{\Pi_x^{(1)}(u)}{\Psi(m(u,T_u))}=\mathcal{P}_{a_1X}^{a_{2}t, \sqrt{2}a_1x}[0,S],
\EQN
with $a_2=\frac{(1+c\varphi)^2}{A_\IF\varphi^2}\left(\alpha_\IF-\frac{c\varphi}{1+c\varphi}\right)$.\\
{\it\underline{Analysis of $\Pi^{(2)}(u)$}}.
It follows from ii) of  Theorem \ref{submain1} that
\BQNY\label{asy1}
\Pi^{(2)}(u)=\pk{\sup_{t\in [0,T_u-S]}X(t)-ct>u}\sim\mathcal{P}_{\mu_\varphi}^{\Omega_\IF(\alpha_\IF-\frac{c\varphi}{1+c\varphi})t}\Psi(m(u,T_u-S)).
\EQNY
By i) \bl{ of  Lemma} \ref{L1} we have that
\BQN\label{ndec1}
\frac{m(u,T_u-S)}{m(u,T_u)}=1+a_u\frac{S}{T_u}(1+o(1)).
\EQN
Then
\BQN\label{asy2}\Pi^{(2)}(u)&\sim& \mathcal{P}_{\mu_\varphi}^{\Omega_\IF(\alpha_\IF-\frac{c\varphi}{1+c\varphi})t}\Psi(m(u,T_u))e^{-\frac{m^2(u,T_u)}{2}\left(\frac{m^2(u,T_u-S)}{m^2(u,T_u)}-1\right)}\nonumber\\
&\sim&\mathcal{P}_{\mu_\varphi}^{\Omega_\IF(\alpha_\IF-\frac{c\varphi}{1+c\varphi})t}\Psi(m(u,T_u))e^{-a_2S}, \quad u\rw\IF.
\EQN
In order to complete the proof of this subcase, we note that combination of (\ref{nin1}), (\ref{new1}) and (\ref{asy2})
leads to
\BQNY
 \mathcal{P}_{a_1X}^{a_{2}t, \sqrt{2}a_1x}[0,S]\leq  \liminf_{u\rw\IF}\frac{\pi_{x,T_u}(u)}{\Psi(m(u,T_u))}\leq\limsup_{u\rw\IF}\frac{\pi_{x,T_u}(u)}{\Psi(m(u,T_u))}\leq \mathcal{P}_{a_1X}^{a_{2}t, \sqrt{2}a_1x}[0,S] +\mathcal{P}_{\mu_\varphi}^{\Omega_\IF(\alpha_\IF-\frac{c\varphi}{1+c\varphi})t}e^{-a_2S}.
\EQNY
 Since  $$\lim_{S\rw\IF}\mathcal{P}_{a_1X}^{a_{2}t, \sqrt{2}a_1x}[0,S]\leq \lim_{S\rw\IF}\mathcal{P}_{a_1X}^{a_{2}t}[0,S]+e^{\sqrt{2}a_1x}<\IF,$$
 which gives the finiteness of the constant,
 then letting $S\rw\IF$ in the above inequalities, we derive
 \BQNY
 \pi_{x,T_u}(u)\sim \mathcal{P}_{a_1X}^{a_{2}t, \sqrt{2}a_1x}\Psi(m(u,T_u)), \quad u\rw\IF.
 \EQNY
{$\diamond$ \underline{Subcase $\varphi\in (0,\IF)$, $\alpha_\IF>1/2$}}.
 Observe that
\BQN\label{nneq2}
&&\pl{X(T_u)-c T_u>u-x, \sup_{0\leq s\leq T_u}\left(X(T_u)-X(s)-c(T_u-s)\right)>u}\nonumber\\
&&\quad \geq \pl{X(T_u)-c T_u>u-x, X(T_u)-c T_u>u}=\Psi(m(u,T_u)).
\EQN
Therefore, in view of (\ref{nneq10}),
\BQN\label{nneq3}\Psi\left(\frac{u-x+cT_u}{\sigma(T_u)}\right)\leq\pi_{x,T_u}(u)\leq  \pi_{0,T_u}(u)+\Psi\left(\frac{u-x+cT_u}{\sigma(T_u)}\right)-\Psi(m(u,T_u)).
\EQN
By (\ref{ratio}) and the fact $\lim_{u\rw\IF}\Omega(u,T_u)=\IF$, in this subcase we have
\BQNY\label{xeq5}
\pi_{0,T_u}(u)\sim \Psi(m(u,T_u))\sim e^{-\frac{x}{A_\IF\varphi^{2\alpha_\IF}}}\Psi\left(\frac{u-x+cT_u}{\sigma(T_u)}\right), \quad u\rw\IF.
\EQNY
which together with (\ref{nneq3}) gives that
$$\pi_{x,T_u}(u)\sim\Psi\left(\frac{u-x+cT_u}{\sigma(T_u)}\right).$$

{$\diamond$ \underline{Case $\varphi=\IF$}}. By the fact that $\lim_{u\rw\IF}\Omega(u,T_u)=\IF$ and (\ref{ratio}), we have
$$ \quad \pi_{0,T_u}(u)\sim \Psi(m(u,T_u))\sim \Psi\left(\frac{u+cT_u-x}{\sigma(T_u)}\right)\quad u\rw\IF.$$
Using the same arguments as given in (\ref{nneq2})-(\ref{nneq3}), we derive
$$\pi_{x,T_u}(u)\sim\Psi\left(\frac{u-x+cT_u}{\sigma(T_u)}\right).$$
This completes the proof. \QED

\subsection{Proof of Theorem \ref{th.0.2}}
 Let $$Y_u(t)=\frac{X(ut)}{u(1+ct)}m(u,t_u), \quad E_2(u)=[t_u/u-\ln m(u,t_u)/m(u,t_u), t_u/u+\ln m(u,t_u)/m(u,t_u)]\cap [0,T_u/u].$$
 Then
  \BQN\label{eq2}
\Pi_3(u)\leq \pi_{0,T_u}(u)\leq\Pi_3(u)+\Pi_4(u),
  \EQN
  where
  \BQN\label{pi34}
  \Pi_3(u)=\pk{\sup_{t\in E_2(u)}Y_u(t)>m(u,t_u)}, \quad \Pi_4(u)=\pk{\sup_{t\in [0,T_u/u]\setminus E_2(u)}Y_u(t)>m(u,t_u)}.
  \EQN
  \gr{In the rest of the proof we shall derive the exact asymptotics of $\Pi_3(u)$.% (by applying Lemma \ref{PPTH0}).
Then we show that $\Pi_4(u)=o(\Pi_3(u))$ as $u\rw\IF$.}
  We distinguish two cases: \kd{$w\in (-\infty,\infty)$ and $w=\infty$.}\\
$\diamond$ \underline{Case $w\in(-\IF,\IF)$ }.\\
 \underline{\it Analysis of  $\Pi_3(u)$}.
\gr{In order to derive the asymptotics of $\Pi_3(u)$, it suffices to check the assumptions of Lemma \ref{PPTH0}. For this, we  observe that}
 $$\Pi_3(u)=\pk{\sup_{t\in [a(u), b(u)]}Y_u\left(\frac{t_u}{u}+\frac{\Delta(u,t_u)t}{u}\right)>m(u,t_u)}, \quad t\in [a(u), b(u)],$$
  with $ a(u)=-\frac{u\ln m(u,t_u)}{\Delta(u,t_u)m(u,t_u)}, \quad b(u)=\min\left(-a(u), \frac{T_u-t_u}{\Delta(u,t_u)}\right)$.
  Moreover, let $ \sigma_u(t)=\sqrt{Var\left(Y_u\left(\frac{t_u}{u}+\frac{\Delta(u,t_u)t}{u}\right)\right)}$ and
 $g(u)=\frac{u^2}{b_u(\Delta(u,t_u))^2}$
 with $b_u$ defined in Lemma \ref{L1}.

In light of Lemma \ref{L1} and \rdd{ Lemma \ref{L2}}, we have
$\sigma_u(0)=1, \ \ 1- \sigma_u(t)\sim \frac{t^2}{g(u)}, t\in [a(u), b(u)],$ and
$$\lim_{u\rw \IF}\sup_{s,t \in [a(u), b(u)], s\neq t}\left|\frac{m^2(u,t_u)\left(1-Corr\left(Y_u\left(\frac{t_u}{u}+\frac{\Delta(u,t_u)t}{u}\right), Y_u\left(\frac{t_u}{u}+\frac{\Delta(u,t_u)t}{u}\right)\right)\right)}{\frac{\sigma^2(\Delta(u,t_u)|s-t|)}{\sigma^2(\Delta(u,t_u))}}-1\right|=0.
$$
\COM{Note that $\frac{\overleftarrow{\sigma}(u^{-1}\sigma^2(u))}{\sigma(u)}\sim \mathbb{Q}u^\beta$, as $u\rw\IF$, with
\BQN\label{beta}\beta=\left\{\begin{array}{cc}
\frac{2\alpha_\IF-1}{\alpha_0}-\alpha_\IF, & \text{if} \ \ \alpha_\IF<1/2\\
-\alpha_\IF, & \text{if} \ \ \alpha_\IF=1/2\\
\frac{2\alpha_\IF-1}{\alpha_\IF}-\alpha_\IF, & \text{if} \ \ \alpha_\IF>1/2.
\end{array}\right.\EQN
Thus $\beta<0$ for all $\alpha_\IF\in (0,1)$ and}
\bl{By ii) of  Lemma \ref{Formula}, we have $$
\lim_{u\rw\IF}\frac{m^2(u,t_u)}{g(u)}=\lim_{u\rw\IF}
\mathbb{Q}\left(\frac{m(u,t_u)\Delta(u,t_u)}{u}\right)^2=0.
$$}
Moreover,
\BQNY
&&\lim_{u\rw\IF}g(u)=\IF,   \quad\lim_{u\rw\IF}\frac{a^2(u)+b^2(u)}{g(u)}=0,
\quad \lim_{u\rw\IF}\frac{a^2(u)+b^2(u)}{g(u)}m(u,t_u)=0,\\
 && y_1=\lim_{u\rw\IF}\frac{a(u)m(u,t_u)}{\sqrt{g(u)}}=-\IF, \quad y_2=\lim_{u\rw\IF}\frac{b(u)m(u,t_u)}{\sqrt{g(u)}}=\lim_{u\rw\IF}\frac{T_u-t_u}{\sqrt{b_u^{-1}A_\IF}t_u^{\alpha_\IF}}
 \frac{u+ct_u}{u}
 =\sqrt{\frac{B}{2AA_\IF}}\frac{(1+ct^*)w}{(t^*)^{\alpha_\IF}}.
\EQNY

Therefore  by i) of Lemma \ref{PPTH0}} we have
\BQNY
\Pi_3(u)&\sim&\mathcal{H}_{\eta_{\alpha_\IF}}\int_{y_1}^{y_2}
e^{-s^2}ds\frac{\sqrt{g(u)}}{m(u,t_u)}\Psi(m(u,t_u))\\
&\sim&\mathcal{H}_{\eta_{\alpha_\IF}}\int_{-\IF}^{\sqrt{\frac{B}{2AA_\IF}}\frac{(1+ct^*)w}{(t^*)^{\alpha_\IF}}}
e^{-s^2}ds\frac{\sqrt{g(u)}}{m(u,t_u)}\Psi(m(u,t_u))\\
&\sim&\mathcal{H}_{\eta_{\alpha_\IF}}\Phi\left(\sqrt{\frac{B}{AA_\IF}}
\frac{(1+ct^*)w}{(t^*)^{\alpha_\IF}}\right)\sqrt{\frac{2A\pi}{B}}\frac{u}{m(u,t_u)\Delta(u,t_u)}\Psi(m(u,t_u)).
\EQNY
\underline{\it Analysis of  $\Pi_4(u)$}. It follows from ii) of Lemma \ref{L1}
that, for $u$ sufficiently large,
  $$\sup_{t\in  [0,T_u/u]\setminus E_2(u)}Var\left(Y_u(t)\right)\leq 1-\frac{B}{4A}m^{-2}(u,t_u)\left(\ln m(u,t_u)\right)^2.$$
  Moreover, by (\ref{reg})
  $$\E{\left(Y_u(t)-Y_u(s)\right)^2}\leq \mathbb{Q}\left(\frac{\sigma^2(u|t-s|)}{\sigma^2(u)}+|t-s|^2\right)\leq \mathbb{Q}|t-s|^\lambda, \ \ s,t\in [0, T_u/u]$$
with $0<\lambda<\min (2\alpha_0, 2\alpha_\IF).$
Applying Piterbarg inequality (Theorem 8.1 in \cite{Pit96}),
we have, as $u\rw\IF$,
\BQN\label{negligible}
\Pi_4(u) \leq\mathbb{Q} (m(u,t_u))^{2/\lambda}\Psi\left(\frac{m(u,t_u)}{\sqrt{1-\mathbb{Q}m^{-2}(u,t_u)(\ln m(u,t_u))^{2}}}\right)=o(\Pi_3(u)).
\EQN
Hence
\BQNY
\pi_{0,T_u}(u)\sim\Pi_3(u), \quad u\rw\IF.
\EQNY

$\diamond$ \underline{Case $w=\IF$}.\\
  \underline{\it Analysis of $\Pi_3(u)$}. Following the same arguments as given for the case $w\in (-\IF,\IF)$, we have that
  $$\Pi_3(u)
\sim\mathcal{H}_{\eta_{\alpha_\IF}}\sqrt{\frac{2A\pi}{B}}\frac{u}{m(u,t_u)\Delta(u,t_u)}\Psi(m(u,t_u)), \quad u\rw\IF.$$
\underline{\it Analysis of  $\Pi_4(u)$}. Observe that
 \BQNY
 \Pi_4(u)
&\leq& \pk{\sup_{t\in [0,M]\setminus E_0(u)}Y_u(t)>m(u,t_u)}+\pk{\sup_{t\in [M,\IF)}Y_u(t)>m(u,t_u)}=:p_1(u)+p_2(u)
 \EQNY
  with $E_0(u)=[t_u-\ln m(u,t_u)/m(u,t_u), t_u+\ln m(u,t_u)/m(u,t_u)
  ]$, $M\in\mathbb{N}$ and $M$ sufficiently large. By the same arguments as given in (\ref{negligible}), we have that
  $p_1(u)=o(\Pi_3(u)),  u\rw\IF.$
  Moreover,
  \BQNY
  p_2(u)\leq \sum_{k=M}^{\IF}\pk{\sup_{t\in [k,k+1]}Y_u(t)
  >m(u,t_u)}.
  \EQNY
  \gr{In order to bound the above sum, we shall apply  Piterbarg inequality  in \cite{Pit96}
for which we
  observe that by Potter's theorem (see, e.g., \cite{BI1989}), we have}
  \BQNY
  \sup_{t\in [k,k+1]} \sqrt{Var\left(Y_u(t)\right)^2}&=&\sup_{t\in [k,k+1]}\frac{\sigma(ut)}{\sigma(t_u)}\frac{1+ct_u}{1+ct}\\
  &\leq&\sup_{t\in [k,k+1]}\mathbb{Q}\left(\frac{t}{t^{*}
  }\right)^{\alpha_\IF+\epsilon}\frac{1+ct_u}{1+ct}\leq \mathbb{Q}k^{\alpha_\IF-1+\epsilon}
  \EQNY
with $0<\epsilon<1-\alpha_\IF$. \gr{Additionally, by (\ref{reg}),} we have for $s,t\in [k,k+1]$ with $k\geq M$,
\BQNY
\E{\left(Y_u(t)-Y_u(s)\right)^2}\leq \mathbb{Q}\left(\frac{\sigma^2(u|t-s|)}{\sigma^2(u)}+|t-s|^2\right)\leq \mathbb{Q}|t-s|^{\lambda}
\EQNY
with $0<\lambda<\min(2\alpha_0, 2\alpha_\IF)$. Thus in light of  Piterbarg inequality  in \cite{Pit96}, we have, for $M$ sufficiently large,
\BQNY
p_2(u)&\leq& \sum_{k=M}^{\IF}\mathbb{Q}\left(m(u,t_u)
\right)^{2/\lambda}\Psi\left(\frac{m(u,t_u)}{\mathbb{Q}k^{\alpha_\IF-1+\epsilon}}\right)\\
&\leq& \mathbb{Q}\left(m(u,t_u)
\right)^{2/\lambda}\Psi\left(\frac{m(u,t_u)}{\mathbb{Q}_1M^{\alpha_\IF-1+\epsilon}}\right)=o(\Pi_3(u)), \quad u\rw\IF.
\EQNY
Consequently, for $M$ sufficiently large,
$\Pi_4(u)=o(\Pi_3(u)), u\rw\IF$. Thus, by (\ref{eq2}),
$$\pi_{0,T_u}(u)\sim\mathcal{H}_{\eta_{\alpha_\IF}}\sqrt{\frac{2A\pi}{B}}\frac{u}{m(u,t_u)\Delta(u,t_u)}\Psi(m(u,t_u)), \quad u\rw\IF.$$
This completes the proof. \QED

\subsection{Proof of Theorem \ref{nth3}} \rdd{Recall that
 \BQN\label{mnth3}
\pi_{x,T_u}(u)&=&\pi_{0,T_u}(u)+\Psi\left(\frac{u-x+cT_u}{\sigma(T_u)}\right)\nonumber\\
&& \ \ -\pl{X(T_u)-c T_u>u-x, \sup_{0\leq s\leq T_u}\left(X(T_u)-X(s)-c(T_u-s)\right)>u}.
\EQN
 The strategy of the proof is the same as used in the proof of Theorem \ref{nth1}, i.e.,
 if $\pi_{0,T_u}(u)$ or $\Psi\left(\frac{u-x+cT_u}{\sigma(T_u)}\right)$ is asymptotically dominating,
 then it determines the  asymptotics of $\pi_{x,T_u}(u)$ as $u\rw\IF$.}
Thus we mostly focus on scenario when this reduction doesn't hold.
 We next provide separate proofs for $\alpha_\IF<1/2$ and $\alpha_\IF\geq 1/2$.\\
{$\diamond$ \underline{Case $\alpha_\IF<1/2$}}. For this case, we distinguish three scenarios.\\
{$\diamond$  \underline{Subcase $\alpha_\IF<1/2$, $\limsup_{u\rw\IF}\frac{T_u-t_u}{\sqrt{u}}<\sqrt{\frac{2A}{B}(1-\alpha_\IF)x}$}}.
We shall prove that
$\pi_{0,T_u}(u)=o\left(\Psi\left(\frac{u-x+cT_u}{\sigma(T_u)}\right)\right)$ as $u\rw\IF$.  \bl{By Theorem \ref{th.0.2} and ii) of  Lemma \ref{Formula},} we have
\BQN\label{xeq1}
\frac{\pi_{0,T_u}(u)}{\Psi\left(\frac{u-x+cT_u}{\sigma(T_u)}\right)}&\sim& \mathbb{Q} u^{-\beta}e^{\frac{1}{2}\left(m^2(u,T_u)-m^2(u,t_u)\right)-\frac{u+cT_u}{\sigma^2(T_u)}x}\nonumber\\
&\sim& \mathbb{Q}e^{\frac{u+cT_u}{\sigma^2(T_u)}\left(\frac{\sigma^2(T_u)}{2(u+cT_u)}\left(m^2(u,T_u)-m^2(u,t_u)\right)-\beta\frac{\sigma^2(T_u)}{u+cT_u}\ln u
-x\right)},
\EQN
with $\beta$ defined in (\ref{beta}).
 \rd{ii) of} Lemma \ref{L1} yields that there exists a sufficiently small $\epsilon>0$ such that
\BQN\label{xeq2}
\frac{\sigma^2(T_u)}{2(u+cT_u)}\left(m^2(u,T_u)-m^2(u,t_u)\right)&=&\frac{\sigma^2(T_u)}{2(u+cT_u)}m^2(u,T_u)\left(1-\frac{m^2(u,t_u)}{m^2(u,T_u)}\right)
\nonumber\\
&\sim&\frac{u}{2(1-\alpha_\IF)}\left(1-\frac{m^2(u,t_u)}{m^2(u,T_u)}\right)\nonumber\\
&\sim& \frac{B}{2A(1-\alpha_\IF)}\left(\frac{T_u-t_u}{\sqrt{u}}\right)^2<x-\epsilon, \quad u\rw\IF.
\EQN
Moreover,
\BQN\label{xeq3}
\frac{\sigma^2(T_u)}{u+cT_u}
\ln u\sim \mathbb{Q}u^{2\alpha_\IF-1}\ln u\rw 0,
\EQN
and
\BQN\label{xeq4}
\frac{u+cT_u}{\sigma^2(T_u)}
\sim \mathbb{Q}u^{1-2\alpha_\IF}\rw\IF.
\EQN
Consequently,
\BQNY
\frac{\pi_{0,T_u}(u)}{\Psi\left(\frac{u-x+cT_u}{\sigma(T_u)}\right)}\leq \mathbb{Q} e^{-\epsilon\frac{u+cT_u}{\sigma^2(T_u)}}\rw 0, \quad u\rw\IF,
\EQNY
which establishes the claim.\\
$\diamond$ \underline{Subcase $\alpha_\IF<1/2$, $\lim_{u\rw\IF}\frac{T_u-t_u}{\sqrt{u}}=\sqrt{\frac{2A}{B}(1-\alpha_\IF)x}$}. Since for this case the above proof doesn't work we shall prove that
$$\pl{X(T_u)-c T_u>u-x, \sup_{0\leq s\leq T_u}\left(X(T_u)-X(s)-c(T_u-s)\right)>u}$$ is negligible compared with $\pi_{0,T_u}(u)+\Psi\left(\frac{u-x+cT_u}{\sigma(T_u)}\right)$, which by (\ref{mnth3}) gives that
$\pi_{x,T_u}(u)\sim \pi_{0,T_u}(u)+\Psi\left(\frac{u-x+cT_u}{\sigma(T_u)}\right)$ as $u\rw\IF$.   We begin with observation that for each $y>0$,
\BQNY
&&\pl{X(T_u)-c T_u>u-x, \sup_{0\leq s\leq T_u}\left(X(T_u)-X(s)-c(T_u-s)\right)>u}\\
&& \ \ = \pl{\overline{X}(T_u)>m(u-x,T_u), \sup_{0\leq s\leq T_u}\left(X(s)-cs\right)>u}\leq \Pi_5(u,y)+\Pi_6(u,y),
\EQNY
where $$\Pi_5(u,y)=\pl{ \sup_{s\in [0, T_u/u]\setminus [t_u/u-y/m(u,t_u),t_u/u+y/m(u,t_u)]}\left(X(us)-cus\right)>u},$$
and
$$\Pi_6(u,y)=\pl{\overline{X}(T_u)>m(u-x,T_u), \sup_{|s-t_u/u|\leq y/m(u,t_u)}\overline{X}(us)>m(u,t_u)}.$$
\underline{ \it Analysis of $\Pi_5(u,y)$}. Observe that
\BQNY\Pi_5(u,y)\leq\pl{ \sup_{s\in [t_u+yu/m(u,t_u), T_u]}\left(X(s)-cs\right)>u}
+\pl{ \sup_{s\in [0, t_u-yu/m(u,t_u)]}\left(X(s)-cs\right)>u}.
\EQNY
Due to the fact that \kd{$\lim_{u\rw\IF}\frac{uy}{m(u,t_u)u^{\alpha_\IF}}=\frac{y\sqrt{A_\IF}(t^*)^{\alpha_\IF}}{1+ct^*}$},
by Theorem \ref{th.0.2}, we have
$$\pl{ \sup_{s\in [0, t_u-yu/m(u,t_u)]}\left(X(s)-cs\right)>u}\sim \Phi(-\sqrt{BA^{-1}}y)\pi_{0,T_u}(u),
\quad u\rw\IF.$$
Following the same arguments as given in the proof of Theorem \ref{th.0.2} for $w=\IF$, we derive that
$$\pl{ \sup_{s\in [t_u+yu/m(u,t_u), T_u]}\left(X(s)-cs\right)>u}\sim \left(1-\Phi(\sqrt{BA^{-1}}y)\right)\pi_{0,T_u}(u), \quad u\rw\IF.$$
Thus
$$\limsup_{u\rw\IF}\frac{\Pi_5(u,y)}{\pi_{0,T_u}(u)}\leq 2\left(1-\Phi(\sqrt{BA^{-1}}y)\right).$$
\underline{\it Analysis of  $\Pi_6(u,y)$}. We have
\BQNY
\Pi_6(u,y)\leq \pl{\sup_{|s-t_u/u|\leq y/m(u,t_u)}\left(\overline{X}(us)+\overline{X}(T_u)\right)>m(u,t_u)+m(u-x, T_u)}.
\EQNY
By  the fact that  $\lim_{u\rw\IF}\frac{T_u-t_u}{\sqrt{u}}=\sqrt{\frac{2A}{B}(1-\alpha_\IF)x}$, we have for $u$ sufficiently large,
$$\inf_{s\in [t_u/u-y/m(u,t_u),t_u/u+y/m(u,t_u)]}|T_u-us|\geq T_u-t_u-yu/m(u,t_u)\geq \sqrt{\frac{A}{B}(1-\alpha_\IF)xu},$$
which together with Lemma \ref{L2} implies that
\BQNY\label{neq1} 2\leq Var\left(\overline{X}(us)+\overline{X}(T_u)\right)&=&4-2(1-\rd{r_u(s, T_u/u)})\nonumber\\
&\leq& 4-\mathbb{Q}\frac{\sigma^2(|T_u-us|)}{\sigma^2(t_u)}\nonumber\\
&\leq& 4-\mathbb{Q}\frac{\sigma^2(|T_u-us|)}{\sigma^2(\Delta(u,t_u))m^2(u,t_u)}\nonumber\\
&\leq& 4-\mathbb{Q}\frac{\sigma^2(\sqrt{u})}{\sigma^2(\Delta(u,t_u))m^2(u,t_u)}
\EQNY
 for all $s\in [t_u/u-y/m(u,t_u),t_u/u+y/m(u,t_u)]$.
Moreover, for  $s,t\in [t_u/u-y/m(u,t_u),t_u/u+y/m(u,t_u)]$ and $u$ large enough
\BQNY1-Corr(\overline{X}(us)+\overline{X}(T_u), \overline{X}(ut)+\overline{X}(T_u))&\leq& \frac{Var(\overline{X}(us)-\overline{X}(ut))}{2\sqrt{Var\left(\overline{X}(us)+\overline{X}(T_u)\right)
Var\left(\overline{X}(ut)+\overline{X}(T_u)\right)}}\\
 &\leq& \frac{1-r_u(s,t)}{2}\leq C\frac{\sigma^2(u|s-t|)}{\sigma^2(ut^*)},
\EQNY
with $C>0$ a fixed constant.
By {\bf AII}, there exists a constant $C_1>0$ such that for $u$ sufficiently large,
\BQNY
\frac{\sigma^2(u|s-t|)}{\sigma^2(ut^*)}=\frac{\sigma^2(\Delta(u,t_u)|u(s-t)/\Delta(u,t_u)|)}{\sigma^2(\Delta(u,t_u))}\frac{\sigma^2(\Delta(u,t_u))}{\sigma^2(ut^*)}\leq C_1\frac{\left|\frac{u}{\Delta(u,t_u)}(s-t)\right|^{\alpha_0/2}}{m^2(u,t_u)}
\EQNY
 for $s,t \in J_k(u):=[t_u/u+k\Delta(u,t_u)/u, t_u/u+(k+1)\Delta(u,t_u)/u]$ with $k\in \mathcal{K}:=\{k: J_k(u)\cap [t_u/u-y/m(u,t_u),t_u+y/m(u,t_u)]\neq \emptyset\}.$
Let $Z_u(t)$ be a family of stationary Gaussian processes with continuous trajectories, unit variances and correlations satisfying
$$Corr(Z_u(t), Z_u(s))=e^{-2CC_1\left(\frac{u}{\Delta(u,t_u)}\right)^{\alpha_0/2}\frac{|t-s|^{\alpha_0/2}}{m^2(u,t_u)}}.$$
 Thus by Slepian's inequality and Lemma 6.1 in \cite{Pit96} we have
\BQNY
\Pi_6(u,y)&\leq& \sum_{k\in \mathcal{K}}\pl{\sup_{s\in J_k(u)}\overline{\overline{X}(us)+\overline{X}(T_u)}>\frac{m(u,t_u)+m(u-x,T_u)}{\sqrt{4-\mathbb{Q}\frac{\sigma^2(\sqrt{u})}{\sigma^2(\Delta(u,t_u))m^2(u,t_u)}}}}\\
&\leq& \sum_{k\in \mathcal{K}}\pl{\sup_{s\in J_k(u)}Z_u(s)>\frac{m(u,t_u)+m(u-x,T_u)}{\sqrt{4-\mathbb{Q}\frac{\sigma^2(\sqrt{u})}{\sigma^2(\Delta(u,t_u))m^2(u,t_u)}}}}\\
&\leq& \mathbb{Q}\frac{u}{m(u,t_u)\Delta(u,t_u)}\pl{\sup_{s\in J_0(u)}Z_u(s)>\frac{m(u,t_u)+m(u-x,T_u)}{\sqrt{4-\mathbb{Q}\frac{\sigma^2(\sqrt{u})}{\sigma^2(\Delta(u,t_u))m^2(u,t_u)}}}}\\
&\leq& \mathbb{Q}\frac{u}{m(u,t_u)\Delta(u,t_u)}\Psi\left(\frac{m(u,t_u)+m(u-x,T_u)}{\sqrt{4-\mathbb{Q}
\frac{\sigma^2(\sqrt{u})}{\sigma^2(\Delta(u,t_u))m^2(u,t_u)}}}\right)(1+o(1)), \quad u\rw\IF,
\EQNY
which implies that for $u$ large enough,
\BQNY
\frac{\Pi_6(u,y)}{\pi_{0,T_u}(u)}&\leq &\mathbb{Q}_1\exp\left({-\frac{(m(u,t_u)+m(u-x,T_u))^2}{8-2\mathbb{Q}
\frac{\sigma^2(\sqrt{u})}{\sigma^2(\Delta(u,t_u))m^2(u,t_u)}}+\frac{m^2(u,t_u)}{2}}\right)\\
&\leq& \mathbb{Q}_1\exp\left({-\frac{(m(u,t_u)+m(u-x,T_u))^2}{8}-\mathbb{Q}_2\frac{\sigma^2(\sqrt{u})}{\sigma^2(\Delta(u,t_u))}+\frac{m^2(u,t_u)}{2}}\right)\\
&\leq& \mathbb{Q}_1\exp\left({\frac{m^2(u,t_u)}{2}\left(1-\frac{1}{4}\left(1+\frac{m(u-x,T_u)}{m(u,t_u)}\right)^2-
2\mathbb{Q}_2\frac{\sigma^2(\sqrt{u})}{\sigma^2(\Delta(u,t_u))m^2(u,t_u)}\right)}\right).
\EQNY
By ii) in Lemma \ref{L1}, we have
\BQNY
1-\frac{1}{4}\left(1+\frac{m(u-x,T_u)}{m(u,t_u)}\right)^2&\sim& 1-\frac{m(u-x,T_u)}{m(u,t_u)}\\
&=&1-\frac{m(u,T_u)}{m(u,t_u)}+\frac{x}{m(u,t_u)\sigma(T_u)}\\
&\sim& -\frac{B}{2A}\left(\frac{T_u-t_u}{u}\right)^2(1+o(1))+\frac{x}{u(1+ct^*)}(1+o(1))=O(u^{-1}), \ \ u\rw\IF.
\EQNY
Moreover,
$$\frac{\sigma^2(\sqrt{u})}{\sigma^2(\Delta(u,t_u))m^2(u,t_u)}\sim \mathbb{Q}_3u^{-\alpha_\IF}. $$
Hence, as $u\rw\IF$,
$$\frac{\Pi_6(u,y)}{\pi_{0,T_u}(u)}\leq \mathbb{Q}_4e^{-\mathbb{Q}_5m^2(u,t_u)u^{-\alpha_\IF}}\sim \mathbb{Q}_4e^{-\mathbb{Q}_5u^{2-3\alpha_\IF}}\rw 0.$$
Consequently,
$$\lim_{y\rw\IF}\limsup_{u\rw\IF}\frac{\Pi_5(u,y)+\Pi_6(u,y)}{\pi_{0,T_u}(u)}=0,$$
implying that $$\pl{\overline{X}(T_u)>m(u-x,T_u), \sup_{0\leq s\leq T_u}\left(X(s)-cs\right)>u}=o(\pi_{0,T_u}(u)).$$
Thus in view of (\ref{mnth3}), we have
\BQNY
\pi_{x,T_u}(u)\sim\pi_{0,T_u}(u)+\Psi\left(\frac{u-x+cT_u}{\sigma(T_u)}\right), \quad u\rw\IF.
\EQNY

$\diamond$ \underline{Subcase $\alpha_\IF<1/2$, $\liminf_{u\rw\IF}\frac{T_u-t_u}{\sqrt{u}}>\sqrt{\frac{2A}{B}(1-\alpha_\IF)x}$}.
 \rdd{First we consider special $T_u$ with $$T_u^*(y)=t_u+\sqrt{\frac{2A}{B}(1-\alpha_\IF)uy}, \quad \text{with} \quad y>x$$ to prove that $\Psi\left(\frac{u-x+cT_u^*(y)}{\sigma(T_u^*(y))}\right)=o\left(\pi_{0,T_u^*(y)}(u)\right)$ as $u\rw\IF$ and then by using the monotonicity of $\Psi\left(\frac{u-x+cT_u}{\sigma(T_u)}\right)$ with respect to $T_u$ we extend this result to all $T_u$ considered.} \\
 Let us first consider the case for $T_u=T_u^*(y)$.
From (\ref{xeq2}),  there exists $\epsilon>0$ such that for $u$ sufficiently large,
\BQNY
\frac{\sigma^2(T_u^*(y))}{2(u+cT_u^*(y))}\left(m^2(u,T_u)-m^2(u,t_u)\right)
&\sim& \frac{B}{2A(1-\alpha_\IF)}u^{-1}(T_u^*(y)-t_u)^2>x+\epsilon
\EQNY
which combined with (\ref{xeq1}), (\ref{xeq3}) and (\ref{xeq4}) leads to, for any $y>x$,
\BQN\label{neq}
\Psi\left(\frac{u-x+cT_u^*(y)}{\sigma(T_u^*(y))}\right)=o\left(\pi_{0,T_u^*(y)}(u)\right), \ \ u\rw \IF.
\EQN
Next we show the monotonicity of  $\Psi\left(\frac{u-x+cT_u}{\sigma(T_u)}\right)$  with respect to $T_u$ sufficiently large.
\rdd{Let
 $t_{u-x}$ denote the maximizer of $m(u-x,\cdot)$ over $(0,\IF)$. Then by ii) in Lemma \ref{L1},  $m(u-x, \cdot)$ is increasing over $[t_{u-x},\IF)$.
 Moreover, one can check that $t_u-t_{u-x}\rw t^*x$ as $u\rw\IF$.  Therefore, $m(u-x,\cdot)$ is increasing over $[t_u+\sqrt{b^{-1}(1-\alpha_\IF)uy},\IF)$ with $y\in (x,\IF)$.\\
  Note that for each $T_u$ satisfying  $\liminf_{u\rw\IF}\frac{T_u-t_u}{\sqrt{u}}>\sqrt{\frac{2A}{B}(1-\alpha_\IF)x}$,  there exists $y\in (x,\IF)$ such that for $u$ large enough
$T_u\geq T_u^*(y)$.} Thus by the monotonicity of $m(u-x,\cdot)$, we have that
\rd{\BQN\label{monoto}
\Psi\left(m(u-x, T_u)\right)\leq \Psi\left(m(u-x,\sigma(T_u^*(y)))\right).
\EQN}
By (\ref{neq})-(\ref{monoto}) and the fact for any $y>0$
$$\pi_{0,T_u^*(y)}(u)\sim \pi_{0,T_u}(u), \quad u\rw\IF,$$
we have
$$\Psi\left(m(u-x,T_u)\right)=o\left(\pi_{0,T_u}(u)\right), \quad u\rw\IF.$$
Hence, $\pi_{x,T_u}(u)\sim \pi_{0,T_u}(u),$ as $ u\rw\IF.$ \\
{$\diamond$ \underline{Case {$\alpha_\IF\geq 1/2$}}}.\\ \rd{ We have
 $$\Psi\left(\frac{u-x+cT_u}{\sigma(T_u)}\right)\sim \Psi\left(\frac{u+cT_u}{\sigma(T_u)}\right)e^{\frac{(u+cT_u)x}{\sigma^2(T_u)}}\sim \mathbb{Q} \Psi\left(\frac{u+cT_u}{\sigma(T_u)}\right), \quad u\rw\IF.$$
Hence by Theorem \ref{th.0.2} we have that
 \BQNY
 \frac{\pi_{0,T_u}(u)}{\Psi(m(u-x,T_u))}\sim\mathbb{Q}_1\frac{\Psi(m(u,t_u))}{ \Psi\left(m(u,T_u)\right)}\frac{u}{m(u,t_u)\Delta(u,t_u)},\quad u\rw\IF.
 \EQNY
 By definition of $t_u$, we have $m(u,t_u)\leq m(u,T_u)$, which  implies that for any $u>0$
 $$\frac{\Psi(m(u,t_u))}{ \Psi\left(m(u,T_u)\right)}\geq 1.$$}
\bl{ Moreover, by  ii) of Lemma \ref{Formula},
 \BQNY
 \frac{u}{m(u,t_u)\Delta(u,t_u)}\sim \mathbb{Q}_2u^{-\beta}\rw\IF, \quad u\rw\IF,
 \EQNY
 where $\beta<0$.}
The above implies that
 $$\Psi(m(u-x,T_u))=o\left(\pi_{0,T_u}(u)\right), \quad u\rw\IF.$$
Hence, $\pi_{x,T_u}(u)\sim \pi_{0,T_u}(u)$, as $u\rw\IF$,
which completes the proof.
 \QED

\def\gtu{g_{u,t}}
\def\xitu{\xi_{u,t}}
\def\chitu{\chi_{u,t}}
\def\zetatu{\zeta_{u,t}}
\def\vtu{\sigma_{u,t}}
\def\bD{E}
\def\sigxiu{\vtu}
\def\cL#1{#1}
\def\H{\mathcal{R}}
\def\coe{\FRAC{C\GAMMA}{1+C\GAMMA}(E)}

{\bf Acknowledgement}:
We would like to thank the referees for their useful comments leading to significant
improvement for the readability of this paper.
K. D\c ebicki
was partially supported by NCN Grant No 2015/17/B/ST1/01102 (2016-2019) whereas
P. Liu was partially supported by  the Swiss National Science Foundation Grant 200021-175752/1.

\bibliographystyle{plain}

 \bibliography{queue}

\end{document}